%
\input amstex
\input amsppt.sty
\magnification=\magstep1
\hsize = 6.5 truein
\vsize = 8 truein
\NoBlackBoxes
\TagsAsMath
\NoRunningHeads

\topmatter
\title Metrics on State Spaces
\endtitle
\author by Marc A. Rieffel
\endauthor
\address\hskip-\parindent
Department of Mathematics \newline University of California
\newline Berkeley, California 94720-3840
\endaddress
\email rieffel\@math.berkeley.edu
\endemail
\date	August 10, 1999
\enddate
\dedicatory This article is dedicated to Richard V\.~Kadison in anticipation
of his completing his seventy-fifth circumnavigation of the sun.
\enddedicatory
\thanks The research reported here was supported in part by National Science
Foundation\newline Grant DMS--96--13833.
\endthanks
\abstract{In contrast to the usual Lipschitz seminorms associated to ordinary
metrics on compact spaces, we show by examples that Lipschitz seminorms on
possibly non-commutative compact spaces are usually not determined by the
restriction of the metric they define on the state space, to the extreme
points of the state space.  We characterize the Lipschitz norms which are
determined by their metric on the whole state space as being those which are
lower semicontinuous.  We show that their domain of Lipschitz elements can be
enlarged so as to form a dual Banach space, which generalizes the situation
for ordinary Lipschitz seminorms.  We give a characterization of the metrics
on state spaces which come from Lipschitz seminorms.  The natural (broader)
setting for these results is provided by the ``function spaces'' of Kadison.
A variety of methods for  constructing Lipschitz seminorms is indicated. }
\endabstract
\endtopmatter

\define\cA{{\Cal A}}
\define\cB{{\Cal B}}
\define\cD{{\Cal D}}
\define\cF{{\Cal F}}
\define\cH{{\Cal H}}
\define\cK{{\Cal K}}
\define\cL{{\Cal L}}
\define\cN{{\Cal N}}
\define\bR{{\Bbb R}}
\define\bC{{\Bbb C}}
\define\a{{\alpha}}
\define\be{{\beta}}
\define\g{{\gamma}}
\define\s{{\sigma}}
\define\e{{\varepsilon}}
\define\lam{{\lambda}}
\define\de{{\delta}}
\let\< = \langle
\let\> = \rangle

\document

In non-commutative geometry (based on $C^*$-algebras), the natural way to
specify a metric is by means of a suitable ``Lipschitz seminorm''.  This idea
was first suggested by Connes \cite{C1} and developed further in \cite{C2,
C3}.  Connes pointed out \cite{C1, C2} that from a Lipschitz seminorm one
obtains in a simple way an ordinary metric on the state space of the
$C^*$-algebra.  This metric generalizes the Monge--Kantorovich metric on
probability measures \cite{KA, Ra, RR}.  In this article we make more
precise the relationship between metrics on the state space and Lipschitz
seminorms.

Let $\rho$ be an ordinary metric on a compact space $X$.  The Lipschitz
seminorm, $L_{\rho}$, determined by $\rho$ is defined on functions $f$ on $X$
by
$$
L_{\rho}(f) = \sup\{|f(x) - f(y)|/\rho(x,y): x \ne y\}  .   \tag 0.1
$$
(This can take value $+\infty$.)  It is known that one can recover $\rho$ from
$L_{\rho}$ by the relationship
$$
\rho(x,y) = \sup\{|f(x)-f(y)|: L_{\rho}(f) \le 1\}.
$$
But a slight extension of this relationship defines a metric, ${\bar \rho}$,
on the space $S(X)$ of probability measures on $X$, by
$$
{\bar \rho}(\mu,\nu) = \sup\{|\mu(f)-\nu(f)|: L_{\rho}(f) \le 1\}. \tag 0.2
$$
This is the Monge--Kantorovich metric.  
The topology which it defines on $S(X)$
coincides with the weak-$*$ topology on $S(X)$ coming from viewing it as the
state space of the $C^*$-algebra $C(X)$.  The extreme points of $S(X)$ are
identified with the points of $X$.  On the extreme points ${\bar \rho}$
coincides with $\rho$.  Thus relationship $(0.1)$ can be viewed as saying that
$L_{\rho}$ can be recovered just from the restriction of its metric ${\bar
\rho}$ on $S(X)$ to the set of extreme points of $S(X)$.

Suppose now that $\cA$ is a unital $C^*$-algebra 
with state space $S(\cA)$, and
let $L$ be a Lipschitz seminorm on $\cA$.  (Precise definitions are given in
Section~2.)  Following Connes \cite{C1, C2}, we define a metric, $\rho$, on
$S(\cA)$ by the evident analogue of $(0.2)$.  We show by simple finite
dimensional examples determined by Dirac operators that $L$ may well not be
determined by the restriction of $\rho$ to the extreme points of $S(\cA)$.

It is then natural to ask whether $L$ is determined by ${\rho}$ on all of
$S(\cA)$, by a formula analogous to $(0.1)$.  
One of our main theorems (Theorem
$4.1$) states that the Lipschitz seminorms for which this is true are exactly
those which are lower semicontinuous in a suitable sense.

For ordinary compact metric spaces $(X,\rho)$ it is known that the space of
Lipschitz functions with a norm coming from the Lipschitz seminorm is the dual
of a certain other Banach space.  
Another of our main theorems (Theorem $5.2$) states that
the same is true in our non-commutative setting, and we give a natural
description of this predual.
We also characterize the metrics on $S(\cA)$ which come from Lipschitz
seminorms (Theorem $9.11$).

We should make precise that we ultimately 
require that our Lipschitz seminorms be
such that the metric on $S(\cA)$ which they determine gives the weak-$*$
topology on $S(\cA)$.  An elementary characterization of exactly when this
happens was given in \cite{Rf}.  (See also \cite{P}.)  This property obviously
holds for finite-dimensional $C^*$-algebras.  It is known to hold in many
situations for commutative $C^*$-algebras, as well as for $C^*$-algebras
obtained by combining commutative ones with finite dimensional ones.  But this
property has not been verified for many examples beyond those.  However in
\cite{Rf} this property was verified for some interesting infinite-dimensional
non-commutative examples, such as the non-commutative tori, and I expect that
eventually it will be found to hold in a wide variety of situations.

Actually, we will see below that the natural setting for our study is the
broader one of order-unit spaces.  The theory of these spaces was launched by
Kadison in his memoire \cite{K1}.  For this reason it is especially
appropriate to dedicate this article to him.  (In \cite{K2} Kadison uses the
terminology ``function systems'', but we will follow \cite{Al} in using the
terminology ``order-unit space'' as being a bit more descriptive of these
objects.)

On the other hand, most of the interesting constructions currently in view of
Lipschitz seminorms on non-commutative $C^*$-algebras, such as those from
Dirac operators, or those in \cite{Rf}, also provide in a natural way
seminorms on all the matrix algebras over the algebras.  
Thus it is likely that
``matrix Lipschitz seminorms'' in analogy with the matrix norms of \cite{Ef}
will eventually be of importance.  
But I have not yet seen how to use them in a
significant way, and so we do not deal with them here.

Let us mention here that a variety of metrics on the state spaces of full
matrix algebras have been employed by the practitioners of quantum mechanics.
A recent representative paper where many references can be found is
\cite{ZS}.  We will later make a few comments relating some of these metrics
to the considerations of the present paper.

The last three sections of this paper will be devoted to a discussion of the
great variety of ways in which Lipschitz seminorms can arise, even for
commutative algebras.  We do discuss here some non-commutative examples, but
most of our examples are commutative.  I hope in a later paper to discuss and
apply some other important classes of non-commutative examples. Some of the
applications which I have in mind will require extending the theory
developed here to quotients and sub-objects.

Finally, we should remark that while we give here considerable attention to
how Dirac operators give metrics on state spaces, 
Connes has shown \cite{C2} that Dirac
operators encode far more than just the metric information.  In particular
they give extensive homological information.  But we do not discuss this
aspect here.

I thank Nik Weaver for suggestions for improvement of the first version
of this article, which are acknowledged more specifically below.

\bigskip
\noindent
{\bf 1.  Recollections on order-unit spaces}

\medskip
We recall \cite{Al} that an order-unit space is a real partially-ordered
vector space, $\cA$, with a distinguished element $e$, the order unit, which
satisfies:

\roster
\item"1)" (Order unit property) For each $a \in \cA$ there is an $r \in \bR$
such that $a \le re$.
\item"2)" (Archimedean property) If $a \in \cA$ and if $a \le re$ for all $r
\in \bR^+$, then $a \le 0$.
\endroster
For any $a \in \cA$ we set
$$
\|a\| = \inf\{r \in \bR^+: -re \le a \le re\}.
$$
We obtain in this way a norm on $\cA$.  In turn, the order can be recovered
from the norm, because $0 \le a \le e$ iff $\|a\| \le 1$ and $\|e-a\| \le 1$.
The primary source of examples consists of the linear 
spaces of all self-adjoint
elements in unital $C^*$-algebras, with the identity element serving as order
unit.  But any linear space of bounded self-adjoint 
operators on a Hilbert space will
be an order-unit space if it contains the identity operator.  We expect that
this broader class of examples will be important for the applications of
metrics on state spaces.

We will not assume that $\cA$ is complete for its norm.  This is important for
us because the domains of Lipschitz norms will be dense, but usually not
closed, in the completion.  (The completion is always again an order-unit
space.)  This also accords with the definition in \cite{Al}.

By a state of an order-unit space $(\cA,e)$ we mean a continuous linear
functional, $\mu$, on $\cA$ such that $\mu(e) = 1 = \|\mu\|$.  States
are automatically positive. We denote the
collection of all the states of $\cA$, i\.e\. the state space of $\cA$, by
$S(\cA)$.  It is a $w^*$-compact convex subset of the Banach space dual,
$\cA'$, of $\cA$.

To each $a \in \cA$ we get a function, ${\hat a}$, on $S(\cA)$ defined by
${\hat a}(\mu) = \mu(a)$.  Then ${\hat a}$ is an affine function on $S(\cA)$
which is continuous by the definition of the $w^*$-topology.
The basic representation theorem of Kadison \cite{K1, K2, K3} (see Theorem
II\.1\.8 of \cite{Al}) says that for any order-unit space the representation
$a \to {\hat a}$ is an isometric order isomorphism of $\cA$ onto a dense
subspace of the space $Af(S(\cA))$ of all continuous affine functions on
$S(\cA)$, equipped with the supremem norm and the usual order on functions
(and with $e$ clearly carried to the constant function $1$).  In
particular, if $\cA$ is complete, then it is isomorphic to all of $Af(S(\cA))$.

Thus we can view the order-unit spaces as exactly the dense subspaces
containing $1$ inside $Af(K)$, where $K$ is any compact convex subset of a
topological vector space.  This provides an effective view from which to see
many of the properties of order-unit spaces.  Most of our theoretical
discussion will be carried out in the setting of order-unit spaces and
$Af(K)$, though our examples will usually involve specific $C^*$-algebras.  We
let $C(K)$ denote the real $C^*$-algebra of all continuous functions on $K$,
in which $Af(K)$ sits as a closed subspace.

It will be important for us to work on the quotient vector space ${\tilde \cA}
= \cA/(\bR e)$.  We let $\|\quad\|^{\sim}$ denote 
the quotient norm on ${\tilde
\cA}$ from $\|\quad\|$.  This quotient norm is easily described.  For $a \in
\cA$ set
$$
\align
\max(a) &= \inf\{r: a \le re\} \\
\min(a) &= \sup\{r: re \le a\},
\endalign
$$
so that $\|a\| = (\max(a))\vee(-\min(a))$.  Then it is easily seen that
$$
\|{\tilde a}\|^{\sim} = (\max(a)-\min(a))/2.
$$

\bigskip
\noindent
{\bf 2.  The radius of the state space}

\medskip
Let $\cA$ be an order-unit space.  Since the term ``Lipschitz seminorm'' has
somewhat wide but imprecise 
usage, we will not use this term for our main objects of precise
study (which we will define in Section 5).  Almost the minimal
requirement for a Lipschitz seminorm is that its null-space be exactly the
scalar multiples of the order unit.  We will use the term ``Lipschitz
seminorm'' in this general sense.  We emphasize that a Lipschitz seminorm will
usually not be continuous for $\|\quad\|$.

Let $L$ be a Lipschitz seminorm on $\cA$.  For $\mu,\nu \in S(\cA)$ we can
define a metric, $\rho_L$, on $S(\cA)$ by
$$
\rho_L(\mu,\nu) = \sup\{|\mu(a)-\nu(a)|: L(a) \le 1\}
$$
(which may be $+\infty$).  Then $\rho_L$ determines a topology on $S(\cA)$.
Eventually we want to require that this topology agrees 
with the weak-$*$ topology.  Since $S(\cA)$ is weak-$*$
compact, $\rho_L$ must then give $S(\cA)$ finite diameter.  We examine this
latter aspect here, in part to establish further notation.

It is actually more convenient for us to work with ``radius'' (half the
diameter), since this will avoid factors of $2$ in various places.  We would
like to use the properties of order-unit spaces to express the radius in terms
of $L$ in a somewhat more precise way than was implicit in \cite{Rf} in its
more general context.  The following considerations \cite{Al} 
will also be used
extensively later.

As in \cite{Rf} and in the previous section, we denote the quotient vector
space $\cA/(\bR e)$ by ${\tilde \cA}$, with its quotient norm
$\|\quad\|^{\sim}$.  But in addition to this norm,
the quotient seminorm ${\tilde L}$ from $L$ is also a
norm on ${\tilde \cA}$, since $L$ takes value $0$ only on $\bR e$.

The dual Banach space to ${\tilde \cA}$ for $\|\quad\|^{\sim}$ is just
${\cA'}^0$, the subspace of $\cA'$ consisting of those $\lam \in \cA'$ such
that $\lam(e) = 0$.  We denote the norm on $\cA'$ dual to  $\|\quad\|$
still by $\|\quad\|$.  The
dual norm on ${\cA'}^0$ is just the restriction of $\|\quad\|$ to ${\cA'}^0$.
If we view $\cA$ as a dense subspace of 
$Af(K) \subseteq C(K)$, then by the Hahn--Banach theorem $\lam$
extends (not uniquely) 
to $C(K)$ with same norm.  There we can take the Jordan decomposition
into disjoint non-negative measures.  Note that for positive measures their
norm on $C(K)$ equals their norm on $\cA$, since $e \in \cA$.  Thus we find
$\mu,\nu \ge 0$ such that $\lam = \mu-\nu$ and $\|\lam\| = \|\mu\| +
\|\nu\|$.  But $0 = \lam(e) = \mu(e) - \nu(e) = \|\mu\| - \|\nu\|$.
Consequently $\|\mu\| = \|\nu\| = \|\lam\|/2$.  Thus if $\|\lam\| \le 2$ we
have $\|\mu\| = \|\nu\| \le 1$.  If $\|\lam\| < 2$ set $t = \|\mu\| < 1$, and
rescale $\mu$ and $\nu$ so that they are in $S(\cA)$.  Then
$$
\lam = t\mu - t\nu = \mu - (t\nu + (1-t)\mu).
$$
Now $(t\nu + (1-t)\mu)$ is no longer disjoint from $\mu$, but we have obtained
the following lemma, which will be used in a number of places.

\proclaim{2.1\ Lemma} The ball $D_2$ of radius $2$ about $0$ in ${\cA'}^0$
coincides with $\{\mu-\nu: \mu,\nu \in S(\cA)\}$.
\endproclaim

Notice that if there is an $a \in \cA$ such that $L(a) = 0$ but 
$a \notin \bR e$, then from this lemma we can find $\mu, \nu \in S(\cA)$
such that $(\mu - \nu)(a) \neq 0$, so that $\rho_L(\mu, \nu) = +\infty$.
Thus our standing assumption that there is no such $a$ serves to reduce the
possibility of having infinite distances. But it does not eliminate this
possibility, as seen by the example of the algebra of smooth 
(or Lipschitz) functions
of compact support on the real line, with constant functions adjoined, and
with the usual Lipschitz seminorm.

\proclaim{2.2\ Proposition} With notation as earlier, the following conditions
are equivalent for an $r \in \bR^+$:

{\rm 1)} For all $\mu,\nu \in S(\cA)$ we have $\rho_L(\mu,\nu) \le 2r$.

{\rm 2)} For all $a \in \cA$ we have $\|{\tilde a}\|^{\sim} \le
rL^{\sim}({\tilde a})$.
\endproclaim

\demo{Proof} Suppose that condition $1$ holds.  Let $a \in \cA$ and $\lam \in
D_2$.  Then by the lemma $\lam = \mu - \nu$ for some $\mu,\nu \in S(\cA)$.
Thus
$$
|\lam(a)| = |(\mu-\nu)(a)| \le L(a)\rho_L(\mu,\nu) \le L(a)2r.
$$
Since $\lam(e) = 0$, thus inequality holds whenever $a$ is replaced by $a +
se$ for $s \in \bR$.  Thus condition $2$ holds.

Conversely, suppose that condition $2$ holds.  Then for any $\mu,\nu \in
S(\cA)$ and $a \in \cA$ with $L(a) \le 1$ we have
$$
|\mu(a)-\nu(a)| = |(\mu-\nu)(a)| \le 2\|{\tilde a}\|^{\sim} \le 2r.
$$
Thus $\rho_L(\mu,\nu) \le 2r$ as desired.\hfill \qed
\enddemo

Of course, we call the smallest $r$ for which the conditions of 
this proposition hold the {\it radius} of $S(\cA)$.

We caution that just because a metric space has radius $r$, it does not follow
that there is a ball of radius $r$ which contains it, as can be seen by
considering equilateral triangles in the plane.
We remark that just because $\rho_L$ gives $S(\cA)$ finite radius, it does not
follow that $\rho_L$ gives the weak-$*$ topology.  Perhaps the simplest
example arises when $\cA$ is infinite dimensional and $L(a) = \|{\tilde
a}\|^{\sim}$.

\bigskip
\noindent
{\bf 3.  Lower semicontinuity for Lipschitz seminorms}

\medskip
Let $L$ be any Lipschitz seminorm on an order-unit space $\cA$. (We
will not at first require that it give $S(\cA)$ finite diameter.)
We would like to show that $L$ and $\rho_L$ contain the same information, and
more specifically that we can recover $L$ from $\rho_L$ as being the usual
Lipschitz seminorm for $\rho_L$.  By this we mean the following.  
Let $\rho$ be
any metric on $S(\cA)$, possibly taking value $+\infty$.  Define $L_{\rho}$
on $C(S(\cA))$ by
$$
L_{\rho}(f) = \sup\{|f(\mu) - f(\nu)|/\rho(\mu,\nu): \mu \ne \nu\}, \tag 3.1
$$
where this may take value $+\infty$.  Let $\text{Lip}_{\rho} = \{f:
L_{\rho}(f) < \infty\}$.  We can restrict $L_{\rho}$ to $Af(S(\cA))$.  In
general, few elements of $Af(S(\cA))$ will be in $\text{Lip}_{\rho}$.
However, on viewing the elements of $\cA$ as elements of $Af(S(\cA))$, 
we have:

\proclaim{3.2\ Lemma} Let $L$ be a Lipschitz seminorm on $\cA$ with corresponding
metric $\rho_L$ on $S(\cA)$.  Then $\cA \subseteq \text{\rm Lip}_{\rho_L}$,
and on $\cA$ we have $L_{\rho_L} \leq L$ , in the sense 
that $L_{\rho_L}(a) \leq L(a)$ for all $a \in \cA$.
\endproclaim

\demo{Proof} For $\mu,\nu \in S(\cA)$ and $a \in \cA$ we have
$$
|{\hat a}(\mu) - {\hat a}(\nu)| = |\mu(a) - \nu(a)| \le L(a)\rho_L(\mu,\nu).
$$
\hfill \qed
\enddemo

For later use we remark that if $L$ and $M$ are Lipschitz seminorms on
$\cA$ and if $M \leq L$, then $\rho_M \geq \rho_L$ in the evident sense.

We would like to recover $L$ on $\cA$ from $\rho_L$ by means of formula
$(3.1)$.  However, the seminorms defined by $(3.1)$ have an important
continuity property:

\definition{3.3\ Definition} Let $\cA$ be a normed vector space, and let $L$
be a seminorm on $\cA$, except that we permit it to take value $+\infty$.
Then $L$ is {\it lower semicontinuous} if for any sequence $\{a_n\}$ in $\cA$
which converges in norm to $a \in \cA$ we have $L(a) \le \liminf\{L(a_n)\}$.
Equivalently, for one, hence every, $t \in \bR$ with $t > 0$, the set
$$
\cL_t = \{a \in \cA: L(a) \le t\}
$$
is norm-closed in $\cA$.
\enddefinition

\proclaim{3.4\ Proposition} Let $\cA$ be an order-unit space, and let $\rho$
be any metric on $S(\cA)$, possibly taking value $+\infty$.
Define $L_{\rho}$
on $C(S(\cA))$ by formula $(3.1)$.  Then $L_{\rho}$ is lower semicontinuous.
Consequently, the restriction of $L_{\rho}$ to any subspace of $C(S(\cA))$,
such as $\cA$ or $Af(S(\cA))$, will be lower semicontinuous.
\endproclaim

\demo{Proof} When we view $L_{\rho}$ as a function of $f$, the formula $(3.1)$
says that $L_{\rho}$ is the pointwise supremum of a collection of functions
(labeled by pairs $\mu,\nu$ with $\mu \ne \nu$) which are clearly continuous
on $C(S(\cA))$ for the supremum norm.  But the pointwise supremum of
continuous functions is lower semicontinuous.\hfill \qed
\enddemo

\example{3.5\ Example} Here is an example of a Lipschitz 
seminorm $L$ whose metric
can be seen to give $S(\cA)$ the weak-$*$ topology,
but which is not lower semicontinuous.  Let $I = [-1,1]$,
and let $\cA = C^1(I)$, the algebra of functions which have continuous first
derivatives on $I$.  Define $L$ on $\cA$ by
$$
L(f) = \|f'\|_{\infty} + |f'(0)|.
$$
For each $n$
let $g_n$ be the function defined by $g_n(t) = n|t|$ for $|t| \le 1/n$, and
$g_n(t) = 1$ elsewhere.  Let $f_n(t) = \int_{-1}^t g_n(s)ds$.  Then the
sequence $\{ f_n\}$
converges uniformly to the function $f$ given by $f(t) = t + 1$.  
But $L(f_n) = 1$
for each $n$, whereas $L(f) = 2$.
\endexample

A substantial supply of examples of lower semicontinuous seminorms can be
obtained from the $W^*$-derivations of Weaver \cite{W2, W3}.  These
derivations will in general have large null spaces, and the seminorms from
them need not give the weak-$*$ topology on the state space.  But many of the
specific examples of $W^*$-derivations which Weaver considers do in fact give
the weak-$*$ topology.  
In terms of Weaver's terminology, which we do not review here, we
have:

\proclaim{3.6\ Proposition} Let $M$ be a von~Neumann algebra and let $E$ be a
normal dual operator $M$-bimodule.  Let $\de: M \to E$ be a
$W^*$-derivation, and denote the domain of $\de$ by $\cL$, so that $\cL$ is
an ultra-weakly dense unital $*$-subalgebra of $M$.  Define a seminorm, $L$,
on $\cL$ by $L(a) = \|\de(a)\|_E$.  Then $L$ is lower semicontinuous, and
$\cL_1 = \{a \in \cL: L(a) \leq 1\}$ is norm-closed in $M$ itself.
\endproclaim

\demo{Proof} Let $\{a_n\}$ be a sequence in $\cL$ which converges in norm to
$b \in M$.  To show that $L$ is lower semicontinuous, it 
suffices to consider the case in which $\{a_n\}$ is contained
in $\cL_1$.  Then the set $\{(a_n,\de(a_n))\}$ 
is a bounded subset of the graph of $\de$ for the norm
$\max\{\|\quad\|_M,\|\quad\|_E\}$.  Since the graph of a $W^*$-derivation is
required to be ultra-weakly closed, and since bounded ultraweakly closed 
subsets are compact for the ultra-weak
topology, there is a subnet which converges ultra-weakly to an element
$(c,\de(c))$ of the graph of $\de$.  Then necessarily $c = b$, so that $b \in
\cL$, and $\de(b)$ is in the ultra-weak closure of $\{\de(a_n)\}$.
Consequently
$L(b) = \|\de(b)\| \le 1$.\hfill \qed
\enddemo

Because of the importance of Dirac operators, it is appropriate to verify
lower semicontinuity for the Lipschitz seminorms which they determine.  This is
close to a special case of Proposition $3.6$, but does not require any kind of
completeness, nor an algebra structure on $\cA$.

\proclaim{3.7\ Proposition} Let $\cA$ be a linear subspace of bounded
self-adjoint operators on a Hilbert space $\cH$, containing the identity
operator.  Let $D$ be an essentially self-adjoint operator on $\cH$ whose
domain, $\cD(D)$, is carried into itself by each element of $\cA$.  Assume
that $[D,a]$ is a bounded operator on $\cD(D)$ for each $a \in \cA$ (so that
$[D,a]$ extends uniquely to a bounded operator on $\cH$).  Define $L$ on $\cA$
by $L(a) = \|[D,a]\|$.  Then $L$ is lower semicontinuous.
\endproclaim

\demo{Proof} Let $\{a_n\}$ be a sequence in $\cA$ which converges in norm to
$a \in \cA$.  Suppose that there is a constant, $k$, such that $L(a_n) \le k$
for all $n$.  For any $\xi,\eta \in \cD(D)$ with $\|\xi\| = 1 = \|\eta\|$
we have
$$
\<[D,a]\xi,\eta\> = \<a\xi,D\eta\> - \<D\xi,a\eta\> 
= \lim\<[D,a_n]\xi,\eta\>.
$$
But $|\<[D,a_n]\xi,\eta\>| \le k$ for each $n$, and so $\|[D,a]\| \le
k$.\hfill \qed
\enddemo

We remark that the Lipschitz seminorms constructed in \cite{Rf} by means of
actions of compact groups are easily seen to be lower semicontinuous.

\bigskip
\noindent
{\bf 4.  Recovering $L$ from $\rho_L$}

\medskip
In this section we show that a lower semicontinuous Lipschitz
seminorm $L$ can be recovered from its metric
$\rho_L$.  But before showing this we would like to emphasize the following
point.  Let $(X,\rho)$ be an ordinary compact metric space, with $\cA$ the
algebra of its Lipschitz functions, with Lipschitz seminorm $L$.  Then $S(\cA)$
consists of the probability measures on $X$, and the points of $X$ correspond
exactly to the extreme points of $S(\cA)$.  The restriction of $\rho_L$
to the extreme points is exactly $\rho$. Thus when one says that one can
recover $L$ from the metric $\rho$, one is saying
that one can recover $L$ from the restriction of $\rho_L$ on $S(\cA)$ to the
extreme points of $S(\cA)$.  However, for the more general situation which we
are considering, it will be false in general 
that we can recover $L$ from the restriction
of $\rho_L$ to the extreme points of $S(\cA)$.  Simple explicit examples will
be given in Section~7.

One of the main theorems of this paper is:

\proclaim{4.1\ Theorem} Let $L$ be a lower semicontinuous Lipschitz
seminorm on an order-unit space
$\cA$, and let $\rho_L$ denote the corresponding metric on $S(\cA)$,
possibly taking value $+\infty$.  Let
$L_{\rho_L}$ be defined by formula $(3.1)$, but restricted to $\cA \subseteq
Af(S(\cA))$.  Then
$$
L_{\rho_L} = L.
$$
\endproclaim

Theorem 4.1 is an immediate consequence of the following theorem, since
we saw that lower semicontinuity coincides with $\cL_1$ being 
norm closed.

\proclaim{4.2\ Theorem} Let $L$ be any Lipschitz seminorm on an
order-unit space $\cA$, and let $\rho_L$ denote the 
corresponding metric on $S(\cA)$.  Let
$L_{\rho_L}$ be defined by formula $(3.1)$, but restricted to $\cA \subseteq
Af(S(\cA))$.  Then $\{a \in \cA: L_{\rho_L}(a) \leq 1\}$ coincides 
with the norm
closure, $\bar \cL_1$, of $\cL_1$ in $\cA$. In particular, $L_{\rho_L}$
is the largest lower semicontinuous seminorm smaller than $L$, and
$\rho_{L_{\rho_L}} = \rho_L$.
\endproclaim

\demo{Proof} (An idea leading to this proof, 
which is simpler than my original 
proof, was suggested to me by Nik Weaver.)
On $\cA'$ we define the seminorm, $L'$ , dual to $L$, by
$$  
      L'(\lambda) = \sup\{|\lambda(a)|: L(a) \leq 1 \}  .
$$
Note that $L'$ takes value $+\infty$ on any $\lambda$ for which
$\lambda(e) \neq 0$, and very possibly on some elements of ${\cA'}^0$
as well. But at any rate we have the following key relationship:

\proclaim{4.3\ Lemma} For $\mu,\nu \in S(\cA)$ we have $\rho_L(\mu,\nu) =
L'(\mu-\nu)$.
\endproclaim

\demo{Proof}
$$
\align
L'(\mu-\nu) &= \sup\{|(\mu-\nu)(a)|: L(a) \le 1\} \\
&= \sup\{|\mu(a)-\nu(a)|: L(a) \le 1\} = \rho_L(\mu,\nu).
\endalign
$$
\hfill \qed
\enddemo

Because $\cL_1$ is already convex and balanced, the bipolar 
theorem \cite{Cw} says that $\bar \cL_1$ is exactly the bipolar
of $\cL_1$. Thus we just need to show 
that $\{a \in \cA: L_{\rho_L}(a) \leq 1\}$ is the bipolar of $\cL_1$.
Now it is clear that the unit $L'$-ball in $\cA'$ is exactly 
the polar \cite{Cw}
of $\cL _1$. This provides the last of the following equivalences.
Let $a \in \cA$. Then:

$L_{\rho_L}(a) \leq 1$ exactly if $|\mu(a) - \nu(a)| \leq \rho_L(\mu, \nu)$ for 
all $\mu, \nu  \in S(\cA)$  , \newline
exactly if $|\lambda(a)| \leq L'(\lambda)$ for 
all $\lambda \in D_2$ (by Lemma 4.3 and Lemma 2.1),
\newline
exactly if $|\lambda(a)| \leq 1$ for 
all $\lambda \in \cA'$ with $L'(\lambda) \leq 1$,   \newline
exactly if $a$ is in the prepolar 
of $\{\lambda: L'(\lambda) \leq 1\}$ (by definition \cite{Cw}),  \newline
exactly if $a$ is in the bipolar of $\cL_1$.

It is clear that $L_{\rho_L}$ is lower semicontinuous, 
that it is the
largest such seminorm smaller than $L$, and that it gives the same metric.
\hfill \qed
\enddemo

Note in particular that if $L$ gives $S(\cA)$ finite diameter, or 
the weak-$*$ topology, then so does $L_{\rho_L}$.

We remark that a sort of dual version of Theorem $4.1$ can be found later in
Theorem $9.7$.

We have the following related considerations.
Suppose again that $L$ is a Lipschitz seminorm 
on an order-unit space $\cA$.  Let ${\bar
\cA}$ denote the completion of $\cA$ for $\|\quad\|$, and let ${\bar \cL}_1$
denote now the closure of $\cL_1$ in ${\bar \cA}$ rather than just 
in $\cA$.  Let
${\bar L}$ denote the corresponding ``Minkowski functional'' on ${\bar \cA}$
obtained by setting, for $b \in {\bar \cA}$,
$$
{\bar L}(b) = \inf\{r \in \bR^+: b \in r{\bar \cL}_1\}.
$$
Since there may be no such $r$, we must allow the value $+\infty$.  With
this understanding,  ${\bar L}$ will be a seminorm on ${\bar \cA}$.  
It is easily seen that ${\bar L}(b) \le 1$ exactly if $b \in
{\bar \cL}_1$, and that ${\bar L}$ is lower semicontinuous because ${\bar
\cL}_1$ is closed.

Up to this point we did not require  lower semicontinuity of $L$.  It's import
is given by:

\proclaim{4.4\ Proposition} Let $L$ be a lower semicontinuous
Lipschitz seminorm on an order-unit space
$\cA$. 
Let ${\bar L}$ on ${\bar \cA}$ be defined as above.  Then ${\bar
L}$ is an extension of $L$, that is, for $a \in \cA$ we have ${\bar L}(a) =
L(a)$. Furthermore, $\rho_{\bar L} = \rho_L$.
\endproclaim

\demo{Proof} Suppose that $a \in \cA$ and $L(a) = 1$.  
Then $a \in \cL_1 \subseteq
{\bar \cL}_1$ and so clearly ${\bar L}(a) \le 1$.  Conversely, if ${\bar L}(a)
\le 1$, then $a \in {\bar \cL}_1$.  Thus there is a sequence $\{a_n\}$ in
$\cL_1$ which converges to $a$, with $L(a_n) \le 1$ for every $n$.  From the
lower semicontinuity of $L$ it follows that $L(a) \le 1$. 
Finally, for $\mu, \nu \in S(\cA)$ we have
$$
\rho_{\bar L}(\mu, \nu) = \sup\{|\mu(a) - \nu(a)|: a \in {\bar \cL_1}\}
= \sup\{|\mu(a) - \nu(a)|: a \in \cL_1\} = \rho_L(\mu, \nu)   .
$$
\hfill \qed
\enddemo

Note in particular that if $L$ gives $S(\cA)$ finite diameter, or 
the weak-$*$ topology, then
so does $\bar L$. However, in general $\bar L$ need not be a Lipschitz
seminorm. For example, let $\cA$ be the algebra of real polynomials
viewed as functions on the interval $[0, 2]$, and let $L$ be the
usual Lipschitz seminorm but defined using only points in $[0,1]$.

\definition{4.5\ Definition} We will call ${\bar L}$ the {\it closure} of
$L$.  We will say that a Lipschitz 
seminorm is {\it closed} if $L = {\bar L}$ (on
the subspace where
${\bar L}$ is finite), or equivalently, if $\cL_1$ is complete for the metric
from $\|\quad\|$.
\enddefinition

Then Proposition 4.4 says that for most purposes we can 
assume that $L$ is closed if convenient.

Suppose now that $L$ is a Lipschitz seminorm 
on $\cA$ which is closed.  On $\cA$ we can
define a new norm, $\left|\|\quad\|\right|$, by
$$
\left|\|a\|\right| = \|a\| + L(a).
$$
It is easily verified that $\cA$ is complete for this new norm.  Suppose that
$\cA$ is a $*$-algebra and $\|\quad\|$ is a $C^*$-norm (this can be
weakened).  Suppose further that $L$ is 
a closed Lipschitz seminorm on $\cA$ which
satisfies the Leibniz inequality.  Then the new norm is a normed-algebra norm,
and so $\cA$ becomes a Banach algebra for the new norm.  In Sections~10 and 11
we will indicate many examples of Lipschitz seminorms satisfying the Leibniz
inequality.  This provides a rich class of examples of Banach algebras which
merit study (even in the cases when they are commutative) along the lines
considered in \cite{BCD, J, W1}.

\bigskip
\noindent
{\bf 5.  The pre-dual of $({\tilde \cA},{\tilde L})$}

\medskip
It has been shown in an increasing variety of situations that the space of
Lipschitz functions with a suitable Lipschitz norm 
is isometrically isomorphic to the
dual of some Banach space.  Some of the history of this phenomenon is sketched
in the notes at the end of chapter~2 of \cite{W1}, or more briefly in
\cite{W2}.  Within the non-commutative setting, Weaver shows in Proposition~2
of \cite{W2} that the domains of $W^*$-derivations (as defined there) are dual
spaces.  However, his $W^*$-derivations can have large null spaces, and they
need not give the weak-$*$ topology on $S(\cA)$.  Nevertheless, Weaver's
approach applies to the non-commutative tori, and gives them the same space
of Lipschitz elements as the approach of the present paper (when combined with
\cite{Rf}).  In fact, Weaver shows in \cite{W3} that for the non-commutative
tori one can also define $\text{Lip}^{\a}$, and 
that $\text{Lip}^{\a}$ is actually
the {\it second} dual of $\text{lip}^{\a}$ when $\a < 1$.

To show within our setting that the space of Lipschitz elements
is the dual of a Banach space, we need to assume that $\rho_L$ gives 
the weak-$*$ topology on $S(\cA)$. As before, 
let $\cL_1 = \{a: L(a) \le 1\}$.  From theorem $1.8$ of
\cite{Rf} we know that $\rho_L$ will give the weak-$*$ topology on $S(\cA)$
exactly if the image of $\cL_1$ in ${\tilde \cA}$ is totally bounded for
$\|\quad\|^{\sim}$.  Equivalently, by theorem $1.9$ of \cite{Rf}, 
$L$ must give $S(\cA)$ finite radius, and for one,
hence all, $t \in \bR$ with $t > 0$, the set
$$
\cB_t = \{a: L(a) \le 1 \text{ and } \|a\| \le t\}
$$
must be totally bounded in $\cA$ for $\|\quad\|$. 
We remark that this implies that if $\{a_n\}$ is a
sequence (or net) in $\cA$ converging pointwise on $S(\cA)$ to $a \in \cA$,
and if there is a constant $k$ such that $\|a_n\| \le k$ and $L(a_n) \le k$
for all $n$, then $a_n$ converges to $a$ in norm.  This is because $\{a_n\}$
is contained in $k\cB_1$ whose closure in the completion ${\bar \cA}$ of $\cA$
is compact.  Let $b$ be any norm limit point of $\{a_n\}$ in ${\bar \cA}$.
Then a subsequence of $\{a_n\}$ converges in norm to $b$. But it still
converges pointwise on $S(\cA)$ to $a$.  Consequently $b = a$, and $a$ is the
only norm limit point of $\{a_n\}$.

We now have in view all the requirements on Lipschitz seminorms
which we need for our present
purposes.  So we now define what we expect is the correct way to
specify metrics on compact non-commutative spaces:

\definition{5.1\ Definition} Let $\cA$ be an order-unit space.  By a {\it
Lip-norm} on $\cA$ we mean a seminorm, $L$, on $\cA$ (taking finite values)
with the following properties:

\roster
\item"1)" For $a \in \cA$ we have $L(a) = 0$ if and only if $a \in \bR e$.
\item"2)" $L$ is lower semicontinuous.
\item"3)" $\{a \in \cA: L(a) \le 1\}$ has image in ${\tilde \cA}$ which is
totally bounded for $\|\quad\|^{\sim}$.
\endroster
\enddefinition

We remark that it is easily checked that the closure (Definition 4.5)
of a Lip-norm is again a Lip-norm.

Within the present setting the fact that the space of Lipschitz elements is a
dual Banach space takes the following form (which requires the Lip-norm
to be closed).

\proclaim{5.2\ Theorem} Let $\cA$ be an order-unit space, and let $L$ be a
Lip-norm on $\cA$ which is closed.  Let $\cK = \{{\tilde a} \in {\tilde \cA}:
{\tilde L}({\tilde a}) \le 1\}$, so that $\cK$ is a compact (convex) set for
$\|\quad\|^{\sim}$.  Then $({\tilde \cA},{\tilde L})$ is naturally
isometrically isomorphic to the dual Banach space of $Af_0(\cK)$, the Banach
space of continuous affine functions on $\cK$ which take value $0$ at $0 \in
{\tilde \cA}$, with the supremum norm.
\endproclaim

\demo{Proof} Let $\cL_1$ and $\cB_t$ be as defined as above.  
Because $L$ is closed, the totally bounded sets $\cB_t$ are complete
for $\|\quad\|$, and so are compact.  From the finite radius considerations of
Section~2 the image of $\cL_1$ in ${\tilde \cA}$ will coincide with the image
of $\cB_t$ for sufficiently large $t$.  Hence the image of $\cL_1$ in ${\tilde
\cA}$ is compact for $\|\quad\|^{\sim}$, not just totally bounded.  But
the image of $\cL_1$ is exactly $\cK$ as defined in the statement of the
theorem.

We can now argue as in the proof of proposition~1 of \cite{W4}.  We include the
argument here in a form specific to our particular situation.

Let $V = Af_0(\cK)$, as defined in the statement of the theorem.  Then from
lemma $4.1$ of \cite{K3} each element of $V$ extends to a linear functional
(not necessarily continuous for $\|\quad\|^{\sim}$) on ${\tilde \cA}$.  But we
still view $V$ as equipped with the uniform norm $\|\quad\|_{\infty}$ from
$C(\cK)$, for which $V$ is complete.  Then for any $f \in V$ we have
$$
\|f\|_{\infty} = \sup\{f({\tilde a}): {\tilde a} \in \cK\} = \sup\{f({\tilde
a}): {\tilde L}({\tilde a}) \le 1\}.
$$
Consequently $\|\quad\|_{\infty}$ is just the dual norm to the norm ${\tilde
L}$ on ${\tilde \cA}$.  But $V$ will usually be much smaller than the entire
dual Banach space of $({\tilde \cA},{\tilde L})$ because of the requirement
that if $f \in V$ then $f$ is continuous on $\cK$.

We let $V'$ denote the dual Banach space to $V$.  We have the evident mapping
$\s$ from ${\tilde \cA}$ to $V'$ defined by $\s({\tilde a})(f) = f({\tilde
a})$.  Use of the Hahn--Banach theorem shows that $Af_0(\cK)$ separates the
points of $\cK$, and from this we see that $\s$ is injective.  Furthermore
$|\s({\tilde a})(f)| = |f({\tilde a})| \le \|f\|_{\infty}{\tilde L}({\tilde
a})$, and so $\|\s\| \le 1$ for the norm ${\tilde L}$ on ${\tilde \cA}$.  In
particular, $\s(\cK) \subseteq (V')_1$, the unit 
ball of $V'$.  From the definitions
of $\s$ and $V$ we see immediately that $\s$ is continuous from $\cK$ to
$(V')_1$ with its weak-$*$ topology from $V$.  Since $\cK$ is compact,
$\s(\cK)$ must be compact for the weak-$*$ topology.  If $\s(\cK)$ were not
all of $(V')_1$, there would be a $\varphi_0 \in (V')_1$ and a weak-$*$
continuous linear functional separating $\varphi_0$ from $\s(\cK)$.  But every
weak-$*$ linear functional comes from $V$.  Thus there would be an $f \in V$
such that
$$
f({\tilde a}) \le 1 < \varphi_0(f)
$$
for every ${\tilde a} \in \cK$.  But the first inequality means that $\|f\|_\infty
\le 1$, and so the second inequality means that $\|\varphi_0\| >1$,
contradicting the assumption that $\varphi_0 \in (V')_1$.  Thus $\s(\cK) =
(V')_1$.  Consequently $\s$ is an isometric isomorphism of $({\tilde
\cA},{\tilde L})$ with $V'$.\hfill \qed
\enddemo

We remark that, if desired, we can make $\cA$ itself into the dual of a Banach
space, in a non-canonical way, as follows.  Let $r$ be the radius of
$(\cA,L)$, and let $\mu$ be any fixed state of $\cA$.  Define an actual norm,
$L_{\mu}$, on $\cA$ by
$$
L_{\mu}(a) = \max\{|\mu(a)|/r,L(a)\}.
$$
Let ${\tilde L}_{\mu}$ be the quotient of $L_{\mu}$ on ${\tilde \cA}$.  It is
clear that ${\tilde L}_{\mu} \ge {\tilde L}$.  But for any given $a \in \cA$
we can find $\a \in \bR$ such that $\|a-\a\| \le r{\tilde L}({\tilde a})$, by
the definition
of radius.  Then
$$
|\mu(a-\a)| \le \|a-\a\| \le r{\tilde L}({\tilde a}),
$$
while $L(a-\a) = {\tilde L}({\tilde a})$.  Consequently ${\tilde
L}_{\mu}({\tilde a}) \le {\tilde L}({\tilde a})$, so that, in fact, ${\tilde
L}_{\mu} = {\tilde L}$.  Thus $(\cA,L_{\mu})$ has $({\tilde \cA},{\tilde L})$
as quotient space.  The quotient map splits by the isometric map ${\tilde a}
\mapsto a - \mu(a)$.  Since $({\tilde \cA},{\tilde L})$ is isometrically
isomorphic to a dual Banach space, it follows easily that $(\cA,L_{\mu})$ is
also.

See also section~2 of \cite{H}, which gives a slightly different approach
because the norm on $\text{Lip}_{\rho}$ is slightly different from that
implicit here.

Let $\cK$ and $V = Af_0(\cK)$ be as in the statement of Theorem $5.2$.  As in
Section~2, the dual of $({\tilde \cA},\|\quad\|^{\sim})$ is ${\cA'}^0$.  By
the finite diameter condition and Proposition $2.2$ each $\lam \in {\cA'}^0$
defines a continuous linear functional on $({\tilde \cA},{\tilde L})$.  Each
such functional is clearly continuous on $\cK$ for its topology from
$\|\quad\|^{\sim}$.  Thus each $\lam \in {\cA'}^0$ defines an element of $V$,
and so we obtain a linear map from ${\cA'}^0$ into $V$.  From Theorem $5.2$
the norm $\|\quad\|_{\infty}$ on $V$ from $C(\cK)$ coincides with the dual
norm $L'$ from $({\tilde \cA},{\tilde L})$. We have the following addition
to Theorem 5.2.

\proclaim{5.3\ Proposition} The image of ${\cA'}^0$ in $Af_0(\cK)$ 
is dense in $Af_0(\cK)$ for its
norm $\|\quad\|_{\infty} = L'$.
\endproclaim

\demo{Proof} Let $\varphi$ be any continuous linear functional on $V$ which is
$0$ on the image of ${\cA'}^0$.  From Theorem $5.2$ every continuous linear
functional on $V$ comes from an element of ${\tilde \cA}$.  If ${\tilde a}$ is
the element of ${\tilde \cA}$ corresponding to $\varphi$, we then have
$\lam({\tilde a}) = 0$ for all $\lam \in {\cA'}^0$, which implies that
${\tilde a} = 0$ so that $\varphi \equiv 0$.  It follows from the Hahn--Banach
theorem that the image of ${\cA'}^0$ is norm dense in $V$.\hfill \qed
\enddemo

\medskip

\bigskip
\noindent
{\bf 6.  Extreme points}

\medskip
Let $L$ be a Lipschitz seminorm on an order-unit 
space $\cA$, and let $\rho_L$ be the
corresponding metric on $S(\cA)$.  Let $E$ denote the set of extreme points of
$S(\cA)$.  Then $E$ need not be a closed subset of $S(\cA)$, but $S(\cA)$ is
the closed convex hull of $E$ by the Krein--Milman theorem.  Of course
$\rho_L$ restricts to a metric on $E$.  We will give explicit examples in the
next section to show that even when $L$ is a Lip-norm 
the restriction of $\rho_L$ to $E$ does not
determine $\rho_L$ or $L$.  Nevertheless, we can try to use the restriction of
$\rho_L$ to define a new Lipschitz seminorm, $L^e$, on $\cA$, by
$$
L^e(a) = \sup\{|\e(a)-\eta(a)|/\rho_L(\e,\eta): \e,\eta \in E,\ \e \ne \eta\}.
$$

\proclaim{6.1\ Proposition} With the above definition, $L^e$ is a lower
semicontinuous Lipschitz seminorm 
on $\cA$, and it is the smallest such on $\cA$ whose
metric on $S(\cA)$ agrees on $E$ with that of $L$. If $L$ is a {\rm Lip}-norm
then so is $L^e$.
\endproclaim

\demo{Proof} From Theorem 4.2 it is clear that we can 
assume that $L$ is lower semicontinuous.
From Theorem $4.1$ we know that any lower semicontinuous Lipschitz 
seminorm, say $L_1$, is
recovered from its metric by a supremum as above, but ranging over all of
$S(\cA)$ rather than just over $E$.  Thus if the metric for $L_1$ agrees with
$\rho_L$ on $E$, we must have $L^e \le L_1$.
By using the argument in the proof of
Proposition $3.4$ it is easily seen that $L^e$ is lower semicontinuous.
Suppose that $L^e(a) = 0$ for some $a \in \cA$.  Recall that $D_2 = \{\lam \in
{\cA'}^0: \|\lam\| \le 2\}$.
\enddemo

\proclaim{6.2\ Lemma} The convex hull of $\{\e - \eta: \e,\eta \in E,\ \e \ne
\eta\}$ is dense in $D_2$ for the weak-$*$ topology.
\endproclaim

\demo{Proof} From Lemma $2.1$ we know that any element of $D_2$ can be
expressed as $\mu-\nu$ for $\mu,\nu \in S(\cA)$.  By the Krein--Milman theorem
each of $\mu,\nu$ can be approximated arbitrary closely in the weak-$*$
topology by convex combinations from $E$, say $\sum \a_j\e_j$ and $\sum
\be_k\eta_k$.  But the difference of such combinations can be expressed as
$$
\sum (\a_j\be_k)(\e_j - \eta_k).
$$
\hfill \qed
\enddemo

From this lemma it is clear that if $L^e(a) = 0$ then $L(a) = 0$, and thus $a
\in \bR e$. Also, it is easy to see that $\rho_{L^e}$ agrees with $\rho_L$
on $E$.

Finally, we must show that if $L$ is a Lip-norm then
the image of $\cK_0 = \{a: L^e(a) \le 1\}$ in ${\tilde \cA}$ 
is totally bounded for $\|\quad\|^{\sim}$.  Notice that this image is
larger than that for $L$, so we can not immediately apply the corresponding
fact for $L$.  Let ${\bar E}$ denote the closure of $E$ in $S(\cA)$.  It is
clear that the supremum defining $L^e$ could just as well be taken over ${\bar
E}$, and so $L^e$ on $\cA$ is just the 
Lipschitz norm for the metric $\rho_L$ restricted
to ${\bar E}$.  Thus $\cK_0$ can be viewed as contained in $\{f \in C({\bar
E}): L^e(f) \le 1\}$, and the latter has totally bounded image in $C({\bar
E})/\bR e$ since it consists of Lipschitz functions for a metric
and ${\bar E}$ is compact.  Thus $\cK_0$ has totally bounded
image in $C({\bar E})/\bR e$.  But the restriction map from $Af(S(\cA))$ to
$C({\bar E})$ is isometric for $\|\quad\|_{\infty}$ since ${\bar E}$ contains
the extreme points.  (See Theorem II\.1\.8 of \cite{Al}.  We are dealing here
with Kadison's smallest separating representation.)  It follows easily that
$\cK_0$ has totally bounded image in ${\tilde \cA}$ as needed.\hfill \qed

\medskip
We remark that if $F$ is any subset of $S(\cA)$ which contains $E$, then we
can use $F$ instead of $E$ to define a Lip-norm $L^F$ just as we defined $L^e$
above.  Then we will have
$$
L^e \le L^F \le L
$$
in the evident sense, with reverse inequalities for the corresponding
metrics.

Suppose that $\cA$ is a dense $*$-subalgebra of a $C^*$-algebra, ${\bar \cA}$,
and that $L$ is a Lip-norm on $\cA$, with corresponding metric $\rho_L$ on
$S(\cA)$.  As above let $E$ denote the set of extreme points of $S(\cA)$.
Assume first that $\cA$ is commutative.  Then $E$ is compact
and ${\bar \cA} \cong C(E)$.  Assume
that $L = L^e$.  Then $L$ is the usual Lipschitz norm coming from the metric
on the compact set $E$ obtained by restricting $\rho_L$ to $E$.  But in this
case we  know that $L$ must then satisfy the Leibniz rule
$$
L(ab) \le L(a)\|b\| + \|a\|L(b).
$$
It is thus natural to ask the general question:

\proclaim{6.3\ Question} What conditions on a {\rm Lip}-norm $L$ on a general
unital $C^*$-algebra imply that $L$ satisfies the Leibniz rule?
\endproclaim

In the next section we will see examples of Lip-norms which do not satisfy $L
= L^e$ and yet satisfy the Leibniz rule.

\bigskip
\noindent
{\bf 7.  Dirac operators and ordinary finite spaces}

\medskip
Connes has shown \cite{C1, C2, C3} that for a compact Riemannian (spin) manifold all
the metric information is contained in the Dirac operator.  This led him to
suggest that for ``non-commutative spaces'', metrics should be specified by
some analogue of Dirac operators.  We explore here some aspects of this
suggestion for finite-dimensional commutative $C^*$-algebras, i\.e\. ordinary
finite spaces.  This will clarify some of the considerations of the previous
sections.  Here and throughout all the rest of this paper, when we say that
an operator $D$ is a ``	Dirac'' operator, this is not meant to indicate any
particular properties of $D$, but rather is meant to indicate how $D$ is
employed, namely to define a Lipschitz seminorm.

Let $X$ be a finite set, and let $\cA = C(X)$.  In order to remain fully in the
setting of the previous sections we take $C(X)$ to consist only of real-valued
functions.  But in the present commutative situation this is not so important
because, unlike the non-commutative case, if one does not know the algebra
structure, the norm for complex-valued functions is still given by a simple
formula in terms of the norm for real-valued functions.  (See e\.g\. lemma~14
of \cite{W2}.)  Consequently we will be a bit careless here about this
distinction.

We will suppose that $\cA$ has been faithfully represented on a
finite-dimensional complex Hilbert space $\cH$.  We suppose given on $\cH$ an
operator $D$ (the ``Dirac'' operator).  It is usual to take $D$ to be
self-adjoint.  But we find it slightly more convenient to take $D$ to be
skew-adjoint.  The two choices are related by a multiplication by $i$, and
give the same metric results.  Following Connes, we define a seminorm, $L$,
on $\cA$ by
$$
L(a) = \|[D,a]\|,
$$
where $[\ ,\ ]$ denotes the usual commutator of operators, and the norm is the
operator norm.  We want $L$ to be a Lip-norm.  Thus we require that if
$[D,a] = 0$ then $a \in \bC I$.  Because we are in a finite-dimensional setting,
$L$ is continuous for $\|\quad\|_{\infty}$, and indeed is a Lip-norm on $\cA$.

From $L$ we obtain a metric, $\rho_L$, on the space $S(\cA)$ of probability
measures on $X$, as well as on its set of extreme points, which is identified
with $X$ itself.  We now give a very simple example to show that $\rho_L$ on
$S(\cA)$ need not agree with the metric obtained from $\rho_L$ on $X$.

\example{7.1\ Example} Consider a three-dimensional commutative $C^*$-algebra,
$\cA$, represented faithfully on a three-dimensional Hilbert space.  Thus we
can identify $\cA$ with the algebra of diagonal matrices in the full matrix
algebra $M_3 = M_3(\bC)$.  We will consider Dirac operators of a special form
which facilitates calculation, namely matrices $D$ in $M_3(\bC)$ of the form
$$
D = \pmatrix 0 &0 &\a \\
0 &0 &\be \\
-\a &-\be &0
\endpmatrix
$$
where $\a > 0$ and $\be > 0$.  We will also restrict to those $f \in \cA$ which are
real, and denote the three values (or diagonal entries) of $f$ by
$(f_1,f_2,f_3)$.  Because $D$ is skew-symmetric, $[D,f]$ is a real symmetric
matrix, whose eigenvalues thus are real.  In fact, we have
$$
[D,f] = \pmatrix 0 &0 &\a(f_3-f_1) \\
0 &0 &\be(f_3-f_2) \\
\a(f_3-f_1) &\be(f_3-f_2) &0
\endpmatrix.
$$
Because of this special form, the eigenvalues are easily calculated, and one
finds that
$$
L(f) = \|[D,f]\| = (\a^2(f_3-f_1)^2 + \be^2(f_3-f_2)^2)^{1/2}.
$$
It is clear from this that if $L(f) = 0$ then $f$ is a constant function.
Thus $L$ defines a Lip-norm on $\cA$.

We now proceed to calculate the corresponding metric on $S(\cA)$.  We first
calculate the dual norm, $L'$, on ${\cA'}^0$, the dual space of ${\tilde \cA}$,
with notation as in the previous sections.  We identify ${\cA'}^0$ with real
diagonal matrices of trace $0$, paired with $\cA$ via the trace.  For $\lam \in
{\cA'}^0$ we denote its components by $\lam = (\lam_1,\lam_2,\lam_3)$.  Of
course
$$
L'(\lam) = \sup\{|\<f,\lam \>|: L(f) \le 1\}.
$$
Now both $|\<f,\lam\>|$ and $L(f)$ are unchanged if we add a constant function
to $f$.  Thus for the supremum defining $L'(\lam)$ we can assume that $f_3 = 0$
always.  Furthermore, we know that $\lam_3 = -(\lam_1+\lam_2)$.  Thus we need
only deal with the first two components of $f$ and $\lam$.  We do this without
changing notation.  Then we see that
$$
L'(\lam) = \sup\{|f_1\lam_1 + f_2\lam_2|: \a^2f_1^2 + \be^2f_2^2 \le 1\}.
$$
But this is just the norm of a functional on a suitable Hilbert space.
Specifically, let $l^2(w)$ be the Hilbert space of functions on a 2-point
space with weight function $w$ given by $(\a^2,\be^2)$.  Then
$$
f_1\lam_1 + f_2\lam_2 = f_1(\lam_1/\a^2)\a^2 + f_2(\lam_2/\be^2)\be^2,
$$
and in this form the norm of the functional is the length of the vector in
$l^2(w)$ defining it.  This gives
$$
\align
L'(\lam) &= ((\lam_1/\a^2)^2\a^2 + (\lam_2/\be^2)^2\be^2)^{1/2} \\
&= (\lam_1^2/\a^2 + \lam_2^2/\be^2)^{1/2}.
\endalign
$$

We now apply this to obtain the metric on $S(\cA)$.  If $\mu,\nu \in S(\cA)$,
then for the evident notation
$$
\rho_L(\mu,\nu) = L'(\mu-\nu) = ((\mu_1-\nu_1)^2/\a^2 +
(\mu_2-\nu_2)^2/\be^2)^{1/2}.
$$

Let $X$ denote the maximal ideal space of $\cA$.  We identify its 3 points
with the 3 extreme points of $S(\cA)$, and label them, corresponding to the
coordinates in $\cA$, by $\de_1,\de_2,\de_3$.  Then from the above formula for
$\rho_L$ we find that the metric on $X$ is given by:
$$
\align
\rho_L(\de_1,\de_2) &= (1/\a^2 + 1/\be^2)^{1/2} \\
\rho_L(\de_1,\de_3) &= 1/\a \\
\rho_L(\de_2,\de_3) &= 1/\be.
\endalign
$$
Define $\g$ by $\rho_L(\de_1,\de_2) = 1/\g$.  Let $L^e$ denote the ordinary
Lipschitz norm on $\cA$ coming from this metric on $X$.  Then
$$
L^e(f) = \max\{|f_1-f_2|\g,|f_1-f_3|\a,|f_2-f_3|\be\}.
$$
Clearly $L^e$ is quite different from $L$.  From Theorem $4.1$ we know that
the metrics on $S(\cA)$ will thus be quite different, even though they agree on
the extreme points.  This is, of course, also easily seen by direct
calculations.
\endexample

We now make some observations in preparation for the next section.  It is
well-known \cite{W1, W2} that the Lipschitz seminorms $L = L_{\rho}$ from ordinary
metrics on a metric space $X$ have a nice relation to the lattice structure of
(real-valued) $C(X)$, namely
$$
L(f\vee g) \le L(f)\vee L(g).
$$

We remark that for the $L$ of the above example this inequality fails.  For
instance, with notation as above, let $f = (1,0,0)$ and $g = (0,1,0)$, so that
$f\vee g = (1,1,0)$.  Then we see that
$$
L(f) = \a,\ L(g) = \be,\ \text{while}\ L(f\vee g) = (\a^2+\be^2)^{1/2}.
$$
(This is related to the counterexample following theorem 16 of \cite{W2}.)

However, it is not difficult to check that the above $L$ does satisfy the
weaker inequality
$$
L(f\vee 0) \le L(f).
$$
In fact, one can prove that this holds for any choice of skew-adjoint $D$ for
the above $\cA$.  To find a counterexample for this weaker inequality one must
take $\cA$ to be $4$-dimensional.  I have not found a systematic way of
constructing a counterexample there, but some examination of what is needed,
followed by some experimentation with MATLAB yields the following (and
related) example:
$$
D = \pmatrix
0 &4 &-1 &0 \\
-4 &0 &2 &-2 \\
1 &-2 &0 &-4 \\
0 &2 &4 &0
\endpmatrix
$$
and $f = (4,2,0,-1)$.

We remark that ordinary Lipschitz norms on compact metric spaces can all be
easily obtained by means of Dirac operators.  I pointed this out in a lecture
in 1993, and the details are indicated after the proof of proposition~8 of
\cite{W2}.  See also the discussion for graphs which we
will give toward the end of Section~11.

\bigskip
\noindent
{\bf 8.  A characterization of ordinary Lipschitz seminorms}

\medskip
Let $X$ be a compact space, let $\rho$ be a metric on $X$ (giving the topology
of $X$), and let $L$ denote the corresponding ordinary Lip-norm on $C(X)$
(permitted to take value $+\infty$).  As just mentioned in the last section,
it is well-known \cite{W1, W2} and easy to prove that $L$ relates nicely to
the lattice structure of $C(X)$ by means of the inequality
$$
L(f\vee g) \le L(f)\vee L(g).
$$
In Weaver's more general setting of domains of $W^*$-derivations he proves
this inequality for $W^*$-derivations of Abelian structure.  (See lemma~12 of
\cite{W2}.)  We show here that the above inequality exactly characterizes the
Lip-norms which are the ordinary Lipschitz seminorms coming from ordinary
metrics on $X$.

We remark that we never assume here that our Lip-norms satisfy the Leibniz
inequality for the algebra structure, namely
$$
L(fg) \le L(f)\|g\| + \|f\|L(g).
$$
But ordinary Lipschitz seminorms do satisfy this inequality. Thus one
consequence of this section is that the above lattice inequality implies the
Leibniz inequality.  On the other hand, the Lip-norm from any ``Dirac''
operator will satisfy the Leibniz inequality, but can easily fail to satisfy
the lattice inequality, as we saw by examples in the previous section.  Thus
the lattice inequality is much stronger than the Leibniz inequality.

However we should point out that for Dirac operators on compact spin
Riemannian manifolds, in spite of their being defined by means of various
partial derivatives and spinors, the corresponding Lip-norms do satisfy the
lattice inequality.  This is because Connes shows \cite{C1, C2, C3} that
the Lip-norms which those Dirac operators define coincide with the ordinary
Lip-norms for the ordinary metrics on the manifolds determined by the
Riemannian metrics.

Recall that for us $C(X)$ consists of real-valued functions.

\proclaim{8.1\ Theorem} Let $X$ be a compact space, let $\cA$ be a dense
subspace of $C(X)$ containing the constant functions, and let $L$ be a
Lip-norm on $\cA$.
Let ${\bar L}$ denote the closure of $L$, viewed as defined
on all of $C(X)$ as in the discussion before Proposition $4.4$,
and thus permitted to take value $+\infty$.  Then the following conditions are
equivalent:

\roster
\item"1." The Lip-norm $L$ is the restriction to $\cA$ of the usual Lipschitz
seminorm corresponding to a metric on $X$ (namely the metric $\rho_L$).

\item"2." For every $f,g \in C(X)$ we have
$$
{\bar L}(f\vee g) \le {\bar L}(f)\vee {\bar L}(g).
$$
\endroster
\endproclaim

The following lemma is somewhat parallel to lemma~13 of \cite{W2}.  For later
use we state it in slightly greater generality than needed immediately.

\proclaim{8.2\ Lemma} Let $\cA$ be a dense subspace of $C(X)$ containing the
constant functions, and closed under the finite lattice operations (i\.e\. if
$f,g \in \cA$ then $f\vee g \in \cA$).  Let $L$ be a {\rm Lip}-norm on $\cA$ which
satisfies the inequality
$$
L(f\vee g) \le L(f)\vee L(g)
$$
for all $f,g \in \cA$.  Let ${\bar L}$ be the closure of $L$, defined on all of
$C(X)$, permitted to take value $+\infty$.  
Let $\cF$ be a bounded subset of $\cA$ for which there is a constant,
$k$, such that $L(f) \le k$ for all $f \in \cF$.  Let $g = \sup\{f
\in \cF\}$.  Then $g \in C(X)$ and ${\bar L}(g) \le k$.
\endproclaim

\demo{Proof} Let $\{g_{\a}\}$ be the net of suprema of finite subsets of
$\cF$.  Then $\{g_{\a}\}$ is contained in $\cA$, and converges up to $f$
pointwise.  By the hypothesis on $L$ we have $L(g_{\a}) \le k$ for all $\a$.
Thus we have
$$
|g_{\a}(x) - g_{\a}(y)| \le k\rho_L(x,y)
$$
for all $\a$ and all $x,y \in X$; that is, $\{g_{\a}\}$ is equicontinuous.  
We can thus
apply the Ascoli theorem \cite{Ru} to conclude that the net $\{g_{\a}\}$ has a
subnet which converges uniformly.  But the limit of this subnet must be $g$,
and so $g$ must be continuous.  Furthermore, from the lower semicontinuity of
${\bar L}$ we must have ${\bar L}(g)\le k$. \qed
\enddemo

\demo{Proof of Theorem 8.1} As indicated above, it is basically well-known,
and not hard to verify, that condition $1$ implies condition $2$.  Suppose
conversely that condition $2$ holds.  For any $x \in X$ let $\rho_L^x$ be the
continuous function on $X$ defined by $\rho_L^x(y) = \rho_L(x,y)$.  Set $S_x =
\{f \in \cA: f(x)= 0,\ L(f) \le 1\}$.  Since $L(f)$ is unchanged when a
constant function is added to $f$, or when $f$ is replaced by $-f$, the
definition of $\rho_L$ can be rewritten as
$$
\rho_L^x(y) = \sup\{f(y): f \in S_x\}.
$$
This means that $\rho_L^x = \sup S_x$.  But $S_x$ is a bounded set in $\cA$ by
the finite radius considerations.  Thus we can apply the above lemma to
conclude that ${\bar L}(\rho_L^x) \le 1$.  Suppose that ${\bar L}(\rho_L^x) =
c < 1$.  Then ${\bar L}((1/c)\rho_L^x) = 1$, and so from the definition of
$\rho_L$ we obtain
$$
(1/c)|\rho_L^x(x) - \rho_L^x(y)| \le \rho_L(x,y),
$$
for all $y \in X$, that is,
$$
\rho_L(x,y) \le c\rho_L(x,y)
$$
for all $y \in X$, which is impossible (unless $X$ has only one point, which
we now do not permit).  Thus ${\bar L}(\rho_L^x) = 1$.

Much as in Section~6, let $L^e$ denote the ordinary Lip-norm on $C(X)$
(permitting value $+\infty$) corresponding to the
restriction of $\rho_L$ as metric on $X$.
(Recall that $X$ is identified with the extreme points of $S(\cA)$.)  As seen
in Proposition~$6.1$, $L^e \le {\bar L}$.  We now show that $L^e = {\bar L}$
because of the inequality in the hypotheses of our theorem (and its
extension in Lemma~$8.2$).  Let $f \in C(X)$, and suppose that $L^e(f) \le
1$.  Thus
$$
|f(x) - f(y)| \le \rho_L(x,y)
$$
for all $x,y \in X$.  In particular
$$
f(x) - \rho_L(x,y) \le f(y).
$$
For each $x \in X$ define $h^x \in C(X)$ by
$$
h^x(y) = f(x) - \rho_L(x,y).
$$
Then the above inequality says that $h^x \le f$ for each $x$.  But it is clear
that $h^x(x) = f(x)$.  Thus $f = \sup\{h^x: x \in X\}$.  Then from the
considerations of the previous paragraph we see that ${\bar L}(h^x) = 1$ for
all $x$.  Thus by Lemma~$8.2$ we have ${\bar L}(f) \le 1$.  It follows that
${\bar L} = L^e$ as desired. \qed
\enddemo

\proclaim{8.3\ Corollary} Let $X$ be a compact space, and let $\cA$ be a dense
subspace of $C(X)$ which contains the constant functions and is closed under
the finite lattice operations.  Let $L$ be a {\rm Lip}-norm on $\cA$, and
suppose that
$$
L(f\vee g) \le L(f)\vee L(g)
$$
for all $f,g \in \cA$.  Then $L$ is the restriction to $\cA$ of the ordinary {\rm
Lip}-norm on $C(X)$ corresponding to the metric $\rho_L$ on $X$.
\endproclaim

\demo{Proof} Let $f,g \in C(X)$.  Then from Lemma $8.2$ we see immediately that
$$
{\bar L}(f\vee g) \le {\bar L}(f)\vee {\bar L}(g).
$$
We can thus apply Theorem $8.1$ to obtain the desired conclusion. \qed
\enddemo

One way of viewing Theorem $8.1$ is that it characterizes the Lip-norms on
commutative $C^*$-algebras which come from the corresponding metric on the
extreme points of $S(\cA)$.  It would be interesting to have a corresponding
characterization for non-commutative $C^*$-algebras, and for general order-unit
spaces.

\bigskip
\noindent
{\bf 9.  Lip-norms from metrics on $S(\cA)$}

\medskip
It is natural to ask which metrics on $S(\cA)$ arise from Lip-norms on $\cA$.
We obtain here a characterization of such metrics.  Many of the steps work for
arbitrary convex sets, and so at first we will work in that setting.  Thus we
let $V$ be any vector space over $\bR$, and we let $K$ be any convex set in
$V$ which spans $V$.  Much as above, let $D_2 = K - K$.  Note that not only is
$D_2$ convex, but it is also balanced, in the sense that if $\lam \in D_2$ and
if $t \in [-1,1]$, then $t\lam \in D_2$.  To see this, note that if $\lam \in
D_2$ then clearly $-\lam \in D_2$, so we only need consider $t \ge 0$.  But
$$
t(\mu-\nu) = \mu - (t\nu + (1-t)\mu),
$$
which is in $D_2$ by the convexity of $K$.  Let $V^0 = \bR D_2$.  Then $V^0$
is a vector subspace of $V$.  In the setting where $K = S(\cA)$ we know that
$V^0$ is a proper subspace of $V$.  Let $M$ be a norm on $V^0$.  Then we can
define a metric, $\rho$, on $K$ by $\rho(\mu,\nu) = M(\mu-\nu)$.  We want to
characterize the metrics which arise in this way.

The most natural property to expect is that $\rho$ be convex (in each
variable), that is:

\definition{9.1\ Definition} We say that a metric $\rho$ on $K$ is {\it
convex} if for every $\mu,\nu_1,\nu_2 \in K$ and $t \in [0,1]$ we have
$$
\rho(\mu,t\nu_1 + (1-t)\nu_2) \le t\rho(\mu,\nu_1) + (1-t)\rho(\mu,\nu_2).
$$
\enddefinition

The metrics coming from norms on $V^0$ are convex because
$$
\mu - (t\nu_1 + (1-t)\nu_2) = t(\mu-\nu_1) + (1-t)(\mu-\nu_2).
$$

Given a metric $\rho$ on $K$, our strategy will be to try to use $\rho$ to
define a norm, $M$, on $V^0$ by first defining it on $D_2$.  Specifically, for
$\lam \in D_2$ we would like to set
$$
M(\lam) = \rho(\mu,\nu)
$$
for $\lam = \mu-\nu$ with $\mu,\nu \in K$.  But we need to know that this is
well-defined.  That is, we need to know that if $\mu,\nu,\mu',\nu' \in K$ and
if $\mu-\nu = \mu'-\nu'$, then $\rho(\mu,\nu) = \rho(\mu',\nu')$.  This can be
rewritten in terms of midpoints so as to appear a bit closer to considerations
of convexity, namely, that if
$$
(\mu + \nu')/2 = (\mu'+\nu)/2 \leqno(9.2)
$$
then $\rho(\mu,\nu) =
\rho(\mu',\nu')$.  This clearly holds for the metrics coming from norms.  One
finds an attractive geometrical interpretation when one draws a picture of
this relation.

\definition{9.3\ Definition} We say that a metric $\rho$ on $K$ is {\it
midpoint-balanced} if whenever equation $(9.2)$ above holds, it follows that
$\rho(\mu,\nu) = \rho(\mu',\nu')$.
\enddefinition

Let us now assume that $\rho$ is midpoint-balanced.  Then $M$ on $D_2$ is
well-defined as above.  We wish to extend it to a norm on $V^0$.  For this to
be possible we first must have the property that if $t \in \bR$, $|t| \le 1$,
and if $\lam \in D_2$, then $M(t\lam) = |t|M(\lam)$.  Now from the definition
of $M$ it is clear that $M(-\lam) = M(\lam)$.  Thus it suffices to treat the
case in which $t \ge 0$.  If $\lam = \mu-\nu$, then
$$
t\lam = t(\mu-\nu) = \mu - (t\nu + (1-t)\mu), 
$$
so that by the definition of $M$ we have $M(t\lam) = \rho(\mu,t\nu + (1-t)\mu)$.
From convexity, $\rho(\mu,t\nu + (1-t)\mu) \le t\rho(\mu,\nu)$.  But also
$t\lam = (t\mu + (1-t)\nu) - \nu$, which gives a similar inequality. Then
from the triangle inequality and convexity we have
$$
\align
   \rho(\mu, \nu) & \le \rho(\mu, t\nu + (1-t)\mu) + \rho(t\nu + (1-t)\mu, \nu) \\
                  & \le t\rho(\mu, \nu) + (1-t)\rho(\mu, \nu) = \rho(\mu, \nu)  .
\endalign
$$

Thus the inequalities must be equalities, and we obtain:

\proclaim{9.4\ Lemma} Let $\rho$ be a metric on $K$ which is convex and
midpoint balanced.  Define $M$ on $D_2$ as above using $\rho$.  Then for any
$\mu,\nu \in S(\cA)$ and $t \in [0,1]$ we have
$$
\rho(\mu,t\nu + (1-t)\mu) = t\rho(\mu,\nu),
$$
and for any $\lam \in D_2$ and $t \in [-1,1]$ we have
$$
M(t\lam) = |t|M(\lam).
$$
\endproclaim

Next, we need that $M$ is subadditive on $D_2$.  This means that if $\lam,\lam'
\in D_2$ and if $\lam + \lam' \in D_2$, then $M(\lam+\lam') \le M(\lam) +
M(\lam)$.  Let $\lam = \mu -\nu$, $\lam' = \mu' - \nu'$.  Then $\lam+\lam' =
(\mu+\mu') -(\nu+\nu')$.  Assuming that $\rho$ is convex and
midpoint-balanced, we obtain from Lemma $9.4$ that
$$
M(\lam+\lam') = 2M((\lam+\lam')/2).
$$
Now $(\lam+\lam')/2 = (\mu+\mu')/2 - (\nu+\nu')/2$, and
$(\mu+\mu')/2,(\nu+\nu')/2 \in S(\cA)$.  Thus
$$
M((\lam+\lam')/2) = \rho((\mu+\mu')/2,(\nu+\nu')/2),
$$
and we see that what we need is:

\definition{9.5\ Definition} We say that a metric $\rho$ on $K$ is {\it
midpoint concave} if for any $\mu,\nu,\mu',\nu' \in K$ we have
$$
\rho((\mu+\mu')/2,(\nu+\nu')/2) \le (1/2)(\rho(\mu,\nu) + \rho(\mu',\nu')).
$$
\enddefinition

Again one finds an attractive geometrical interpretation when one draws a
picture of this inequality.  From the discussion above we now know that:

\proclaim{9.6\ Lemma} Let $\rho$ be a metric on $K$ which is convex, midpoint
balanced, and midpoint concave.  Define $M$ on $K$ as above.  If $\lam,\lam'
\in D_2$ and if $\lam+\lam' \in D_2$, then
$$
M(\lam+\lam') \le M(\lam) + M(\lam').
$$
\endproclaim

\proclaim{9.7\ Theorem} Let $\rho$ be a metric on the convex subset $K$ of
$V$, and let $V^0 = \bR D_2 = \bR(K-K)$.  Then there is a norm, $M$, on $V^0$
such that $\rho(\mu,\nu) = M(\mu-\nu)$ for all $\mu,\nu \in K$, if and only if
$\rho$ is convex, midpoint balanced, and midpoint concave.  The norm $M$ is
unique.
\endproclaim

\demo{Proof} The uniqueness is clear since $V^0 = \bR(K-K)$.  We have seen
above that the conditions on $\rho$ are necessary.  We now show that they are
sufficient.  We let $M$ be defined on $D_2 = K-K$ as above.  For any $\lam \in
V^0$ there is a $t > 0$ such that $t\lam \in D_2$.  We want to extend $M$ to
$V^0$ by setting
$$
M(\lam) = t^{-1}M(t\lam).
$$
From Lemma $9.4$ it is easily seen that $M$ is well-defined, and furthermore
that $M(s\lam) = |s|M(\lam)$ for all $s \in \bR$ and $\lam \in V^0$. 
The subadditivity of $M$ then follows easily from Lemma 9.6.   \qed
\enddemo

We now want to apply the above ideas to $S(\cA)$ for an order-unit space
$\cA$.  Note that the $V^0$ of just above is then the ${\cA'}^0$ of earlier.
We will need the following theorem, which does not involve the above ideas.

\proclaim{9.8\ Theorem} Let $\cA$ be an order-unit space, and let $M$ be a
norm on ${\cA'}^0$.  Define a metric, $\rho$, on $S(\cA)$ by
$$
\rho(\mu,\nu) = M(\mu-\nu).
$$
If the $\rho$-topology coincides with the weak-$*$ topology on $S(\cA)$, then
$$
M = (L_{\rho})'
$$
on ${\cA'}^0$.
\endproclaim

\demo{Proof} Since $\text{Lip}_\rho$ is a subspace of
$C(S(\cA))$, we can set $\cA_L = (\text{Lip}_\rho) \cap Af(S(\cA))$.  Note that
$\cA_L$ need not be contained in $\cA$ unless $\cA$ is complete.  Initially it
is not clear how big $\cA_L$ is.  Parallel to our earlier notation, let $V$
denote the normed space ${\cA'}^0$ with norm $M$.  Note that $V$ need not be
complete.  Let $V'$ denote the Banach space dual of $V$, with dual norm $M'$.
Fix any $\nu_0 \in S(\cA)$.  For any $\varphi \in V'$ define a function,
$\tau(\varphi)$, on $S(\cA)$ by
$$
\tau(\varphi)(\mu) = \varphi(\mu-\nu_0).
$$
Then for $\mu,\nu \in S(\cA)$ we have
$$
|\tau(\varphi)(\mu) - \tau(\varphi)(\nu)| = |\varphi(\mu-\nu)| \le
M'(\varphi)M(\mu-\nu) = M'(\varphi)\rho(\mu,\nu).
$$
Thus $\tau(\varphi) \in \text{Lip}_{\rho}$ and $L_{\rho}(\tau(\varphi)) \le
M'(\varphi)$.  In particular, $\tau(\varphi)$ is continuous on $S(\cA)$ since
$\rho$ gives the weak-$*$ topology.  Furthermore it is easily seen that
$\tau(\varphi)$ is affine on $S(\cA)$.  Thus $\tau(\varphi) \in \cA_L$.
Consequently $\tau$ is a norm-non-increasing linear map from $(V',M')$ to
$(\cA_L,L_{\rho})$.  Let ${\tilde \tau}$ denote $\tau$ composed with the map
from $\cA_L$ to ${\tilde \cA}_L$.  Then it is easily seen that ${\tilde \tau}$
does not depend on the choice of $\nu_0$. We now need:

\proclaim{9.9\ Lemma} Let ${\bar \cA} = Af(S(\cA))$, the completion of $\cA$
for $\|\quad\|$, so that $\cA_L \subseteq {\bar \cA}$.  Then $\cA_L$ is dense
in ${\bar \cA}$.
\endproclaim

\demo{Proof} Since $\bR e \subseteq \cA_L$, it suffices to show that ${\tilde
\cA}_L$ is dense in ${\bar \cA}^{\sim}$.  Let $\lam \in D_2 \subseteq {\cA'}^0
= ({\bar \cA}^{\sim})'$.  Suppose that $\lam(\cA_L) = 0$.  Let $\lam= \mu-\nu$
with $\mu,\nu \in S(\cA)$.  For any $\varphi \in V'$ we have $\tau(\varphi)
\in \cA_L$, so
$$
0 = \lam(\tau(\varphi)) = \mu(\tau(\varphi)) -
\nu(\tau(\varphi)) 
= \varphi(\mu-\nu_0) - \varphi(\nu-\nu_0) = \varphi(\lam)  .
$$
Since this is true for all $\varphi \in V'$, it follows that $\lam = 0$.
Since $D_2$ spans ${\cA'}^0$, an application of the Hahn--Banach theorem now
shows that $\cA_L$ is dense on ${\bar \cA}$. \qed
\enddemo

Now let $f \in \cA_L$.  We seek to define a linear functional, $\s(f)$, on
${\cA'}^0$ related to the $\s$ in the proof of Theorem $5.2$.  We first try to
define $\s$ on $D_2$ by
$$
\s(f)(\lam) = f(\mu) - f(\nu),
$$
where $\lam = \mu-\nu$ for $\mu,\nu \in S(\cA)$.  But we need to show that
$\s(f)$ is well-defined.  We argue much as we did before Definition $9.3$.  If
also $\lam = \mu_1 - \nu_1$ for $\mu_1,\nu_1 \in S(\cA)$, then $(\mu+\nu_1)/2
= (\mu_1+\nu)/2$.  But these are elements of $S(\cA)$ and so
$$
f((\mu+\nu_1)/2) = f((\mu_1+\nu)/2).
$$
But from the fact that $f$ is affine it now follows that
$$
f(\mu) - f(\nu) = f(\mu_1) - f(\nu_1).
$$
Thus $\s(f)$ is well-defined on $D_2$.  We now need to know that $\s(f)$ is
``linear'' on $D_2$.  The proof that $\s(f)(t\lam) = t\s(f)(\lam)$ for $t \in
[-1,1]$ is similar to the proof of Lemma $9.4$.  The proof that
$\s(f)(\lam+\lam_1) = \s(f)(\lam) + \s(f)(\lam_1)$ if $\lam+\lam_1 \in D_2$ is
similar to the argument just before Definition $9.5$.  The proof that $\s(f)$
then extends to a linear functional on ${\cA'}^0$ is similar to the 
arguments in the proof of
Theorem $9.7$.  For $\lam = \mu -\nu$ with $\mu,\nu \in S(\cA)$ we have
$$
|\s(f)(\lam)| = |f(\mu)-f(\nu)| \le L_{\rho}(f)\rho(\mu,\nu) =
L_{\rho}(f)M(\mu-\nu) = L_{\rho}(f)M(\lam).
$$
It follows that $\s(f) \in V'$ and $M'(\s(f)) \le L_{\rho}(f)$.  Thus $\s$ is
a norm-non-increasing linear map from $(\cA_L,L_{\rho})$ to $(V',M')$.  Note
that the constant functions are in the kernel of $\s$, so that $\s$ determines
a norm-non-increasing linear map from $({\tilde \cA}_L,{\tilde L}_{\rho})$ to
$(V',M')$.  But for $f \in \cA_L$ we have
$$
\tau(\s(f))(\mu) = \s(f)(\mu-\nu_0) = f(\mu)-f(\mu_0).
$$
Consequently ${\tilde \tau}({\tilde \s}({\tilde f})) = {\tilde f}$.
Similarly, for $\varphi \in V'$ and $\lam = \mu-\nu$ we have
$$
{\tilde \s}({\tilde \tau}(\varphi))(\lam) = \tau(\varphi)(\mu) -
\tau(\varphi)(\nu) 
= \varphi(\mu-\nu_0) - \varphi(\nu-\nu_0) = \varphi(\lam),
$$
so that ${\tilde \s}({\tilde \tau}(\varphi)) = \varphi$.  Thus ${\tilde \s}$
and ${\tilde \tau}$ are inverses of each other.  Since they are 
norm-non-increasing, we obtain:

\proclaim{9.10\ Lemma} The map ${\tilde \tau}$ is an isometric isomorphism of
$(V',M')$ onto $(\cA_L,L_{\rho})$, with inverse ${\tilde \s}$.
\endproclaim

We can now complete the proof of Theorem $9.8$. Since $\cA_L$ is 
dense in ${\bar \cA}$ by Lemma $9.9$, for any
$\lam \in V'$ we have
$$
(L_{\rho})'(\lam) = 
\sup\{\lam({\tilde \tau}(\varphi)): L_{\rho}({\tilde
\tau}(\varphi)) \le 1\} 
= \sup\{\varphi(\lam): M'(\varphi) \le 1\} = M(\lam)  .
$$
\qed
\enddemo

Putting together the various pieces of this section, we obtain:

\proclaim{9.11\ Theorem} Let $\cA$ be an order-unit space, and let $\rho$ be a
metric on $S(\cA)$ which gives the weak-$*$ topology.  Then $\rho$ comes from
a {\rm Lip}-norm $L$ on $\cA$ via the relation
$$
\rho(\mu,\nu) = L'(\mu-\nu)
$$
if and only if $\rho$ is convex, midpoint balanced, and midpoint convex.
\endproclaim

Nik Weaver has suggested to me the following alternative treatment of the
material of this section. Let $V$, $K$, and $V^0$ be as at the 
beginning of this section. 

\definition{9.12\ Definition} We say that a metric $\rho$ on $K$ is
{\it linear} if for every $\mu, \nu \in K$, every $v \in V^0$, and every
$t \in \bR^+$ such that $\mu + tv$ and $\nu + v$ are in $K$ we have
$$
\rho(\mu, \mu + tv) = t\rho(\nu, \nu+v)  .
$$
\enddefinition

It is easily seen that if $\rho$ comes from a norm on $V^0$ then
$\rho$ is linear. Conversely, if $\rho$ is linear, define a norm, $M$,
on $V^0$ by
$$
M(v) = \rho(\mu, \mu+tv)/t
$$
for any $\mu \in K$ and any $t \in \bR^+$ such that $\mu + tv \in K$.
One checks that $M$ is well-defined and is indeed a norm. Furthermore,
$\rho$ comes from $M$.

Weaver also points out that if $V$ is a locally convex topological
vector space and if $K$ is compact, then for a suitable definition
of $\rho$ being compatible with the topology, one can show that
when $\rho$ is linear and compatible, then $K$ is isometrically
isomorphic to $S(Af(K))$ when the latter is given the metric
coming from the Lipschitz seminorm on $Af(K)$ coming from $\rho$.

It is not clear that examples will come up where it is actually
useful to apply the considerations of this section in order to
obtain Lip-norms. Until such examples arise, it will not be
clear whether my version or Weaver's will be the more useful.

\bigskip
\noindent
{\bf 10.  Musings on metrics}

\medskip
Since the theory in the previous sections worked for order-unit spaces, which
need not be algebras, the Leibniz inequality played no significant role there.
Indeed, even when one has an algebra, I have not seen how to make effective
use of the Leibniz inequality.  Nevertheless, most constructions of Lipschitz
seminorms which I have seen in the literature seem to provide ones which do
satisfy the Leibniz inequality.  We will briefly explore here a variety of
such constructions, and the relationships between them.  
Our interest will be on
seeing general patterns, and we will not try to deal carefully with the many
technical issues which arise.  Thus we will be less precise than in the
previous sections.

A very natural way to look for Lipschitz seminorms, closely related to Weaver's
$W^*$-derivations \cite{W2}, goes as follows.  Let $\cA$ be a unital algebra
and let $(\Omega,d)$ be a first-order differential calculus for $\cA$.  Thus
$\Omega$ (which is also often denoted $\Omega^1$) is an $\cA$-$\cA$-bimodule,
and $d$ is an $\Omega$-valued derivation on $\cA$, that is, a linear map from
$\cA$ into $\Omega$ which satisfies the Leibniz identity
$$
d(ab) = (da)b + a(db).
$$
We do not require that the range of $d$ generates $\Omega$.  Suppose now that
$\cA$ is in fact a normed algebra, and that we have a bimodule norm, $N$, on
$\Omega$ (for the norm $\|\quad\|$ on $\cA$), that is, a norm such that
$$
N(a\omega b) \le \|a\|N(\omega)\|b\|
$$
for $a,b \in \cA$ and $\omega \in \Omega$.  Define a seminorm $L$ on $\Omega$
by
$$
L(a) = N(da).
$$
It is easily seen that $L$ satisfies the Leibniz inequality.  Since $d1 = 0$,
we have $L(1) = 0$.  Of course, without further hypotheses the null-space of
$L$ may be much bigger.  (We should mention that not all seminorms satisfying
the Leibniz inequality can be constructed in this way---see the discussion in
\cite{BC}.)

There is a universal first-order differential calculus for any unital
algebra $\cA$ \cite{Ar,
C2}.  We approach this in a way which emphasizes more than usual those
differential calculi which are {\it inner}, since at least conceptually that
is what Dirac operators give, as we will see shortly.  We form the algebraic
tensor product
$$
\Omega_1^u = \cA \otimes \cA,
$$
with bimodule structure defined as usual by $a(b \otimes c)d = ab \otimes
cd$.  We define $d$ by
$$
da = 1 \otimes a - a \otimes 1.
$$

\definition{10.1\ Definition} A first-order calculus $(\Omega,d)$ is {\it
inner} if there is a $\omega_0 \in \Omega$ such that
$$
da = \omega_0a - a\omega_0.
$$
\enddefinition

Then the calculus $(\Omega_1^u, d)$ defined above is inner, with $\omega_0 =
1 \otimes 1$.  Note that here $\omega_0$ may not be in the sub-bimodule
generated by the range of $d$.  This is an indication of why we do not
require this generation property.  It is simple to verify:

\proclaim{10.2\ Proposition} The inner first-order calculus $(\Omega_1^u,d,1
\otimes 1)$ is universal among inner first-order differential calculi over
$\cA$, in the sense that if $(\Omega',d',\omega'_0)$ is any other inner
first-order differential calculus, then there is a bimodule homomorphism $\Phi:
\Omega_1^u \rightarrow \Omega'$ such that $\Phi(da) = d'a$ and $\Phi(1 \otimes
1) = \omega'_0$.  In particular,
$$
\Phi(a \otimes b) = a\omega'_0 b
$$
for $a,b \in \cA$.  If $\Omega'$ is generated by $\omega'_0$ as bimodule,
then $\Phi$ is surjective, so that $\Omega'$ is a quotient of $\Omega_1^u$.
\endproclaim

\proclaim{10.3\ Proposition} Any first-order differential calculus is
contained in an inner first-order calculus.
\endproclaim

\demo{Proof} Let $(\Omega,d)$ be a first-order calculus.  Set ${\bar \Omega}
= \Omega \oplus \cA$ as left $\cA$-module, set ${\bar d}a = da \oplus 0$, and
set ${\bar \omega}_0 = 0 \oplus 1$.  We must extend the right action of $\cA$
on $\Omega$ to a right action on ${\bar \Omega}$ such that ${\bar d}a = {\bar
\omega}_0a - a{\bar \omega}_0$.  Thus it
is clear that we must set $(0 \oplus 1)a = {\bar \omega}_0a = da \oplus 0 +
a{\bar \omega}_0 = da \oplus a$, and so
$$
(\omega,b)a = (\omega a + bda,ba).
$$
It is simple to check that this gives the desired structure. \qed
\enddemo

Now let $\Omega^u$ denote the sub-bimodule of $\Omega_1^u$ generated by the
range of $d$, and so spanned by elements of the form
$$
adb = a \otimes b - ab\otimes 1.
$$
Let $(\Omega',d')$ be a first-order differential calculus which is not
inner.  Expand it to an inner calculus by the construction of the previous
proposition, and then restrict $\Phi$ of that proposition to $\Omega^u$.  It is
clear from the construction that $\Phi$ will carry $\Omega^u$ into $\Omega'$,
where $\Omega'$ is viewed as a sub-bimodule of its expansion.  We obtain in
this way:

\proclaim{10.4\ Proposition} The calculus $(\Omega^u,d)$ is universal among
all first-order differential calculi over $\cA$, in the sense that if
$(\Omega',d')$ is any other first-order differential calculus, then there is a
bimodule homomorphism $\Phi: \Omega^u\rightarrow \Omega$ such that $\Phi(da) =
d'a$.  If $\Omega'$ is generated by the range of $d'$ as bimodule, then
$\Phi$ is surjective, so that $\Omega'$ is a quotient of $\Omega^u$.
\endproclaim

We notice that if $(\Omega,d)$ is any first-order differential calculus and if
$\cN$ is any sub-bimodule of $\Omega$, then we obtain a calculus
$(\Omega/\cN,d')$ where $d'$ is the composition of $d$ with the canonical
projection of $\Omega$ onto $\Omega/\cN$.  However, unlike the universal
calculus, there may now be many more elements $a$ for which $da = 0$ beyond
the scalar multiples of $1$.

Let us examine briefly what the above looks like when $\cA = C(X)$ for a
compact space $X$.  Then $\Omega_1^u(= \cA \otimes \cA)$ is naturally viewed
as a dense sub-bimodule, in fact subalgebra, of $C(X \times X)$.  The bimodule
actions are, of course,
$$
(fF)(x,y) = f(x)F(x,y),\ (Ff)(x,y) = F(x,y)f(y),
$$
and $\omega_0 = 1 \otimes 1$ is the constant function $1$, so that $d$ is
given by
$$
(df)(x,y) = f(y) - f(x).
$$
Then $\Omega^u$ is spanned by the $fdg$, where
$$
(fdg)(x,y) = f(x)(g(y) - g(x)).
$$
Thus the elements of $\Omega^u$ take value $0$ on the diagonal, $\Delta$, of
$X \times X$, and consequently $\Omega^u \subseteq C_{\infty}(X \times
X\smallsetminus\Delta)$.  In fact it is easy to see that $\Omega^u$ is a dense
subalgebra of $C_{\infty}(X \times X\smallsetminus\Delta)$.

Let $\rho$ be an ordinary metric on $X$ (giving the topology of $X$).  View
$\rho$ as a strictly positive function on $X \times X\smallsetminus\Delta$,
and let $\g = \rho^{-1}$.  Then $\g$ is a continuous function on $X \times
X\smallsetminus\Delta$, but $\g$ is unbounded if $X$ is not finite.  Let $C(X
\times X\smallsetminus\Delta)$ denote the algebra of continuous possibly-unbounded
functions on $X \times X\smallsetminus\Delta$.  Then $C(X \times
X\smallsetminus\Delta)$ can be viewed as the algebra of operators affiliated with
the $C^*$-algebra $C_{\infty}(X\times X\smallsetminus\Delta)$ in the sense
studied by Baaj \cite{Ba} and Woronowicz \cite{Wo}.  In an evident way $C(X
\times X\smallsetminus\Delta)$ is an $\cA$-$\cA$-bimodule, containing $\g$.

There are now two routes which we can take.  One is to consider the
inner-derivation, $d_{\g}$, defined by $\g$.  Thus
$$
(d_{\g}f)(x,y) = \g(x,y)f(y) - f(x)\g(x,y) = (f(y)-f(x))/\rho(x,y).
$$
Then we can consider bimodule norms, possibly taking value $+\infty$, on $C(X
\times X\smallsetminus\Delta)$, as a way to obtain Lipschitz norms on $\cA$.
The other route is to use $\g$ (or $\rho$) to directly define norms on
$C_{\infty}(X \times X\smallsetminus\Delta)$.  For the first route the most
obvious norm is the supremum norm, which leads to the usual definition of the
Lipschitz seminorm for a metric space.

However, we choose to explore further the second route.  (But most of what we
find will have a fairly evident reinterpretation in terms of the first
route.)  There is a large variety of ways to obtain bimodule norms on
$C_{\infty}(X \times X\smallsetminus\Delta)$.  The one which gives the usual
definition of the Lipschitz seminorm for a metric is clearly
$$
N(F) = \|\g F\|_{\infty},
$$
permitted to take value $+\infty$.  But here are some others.  Let $m$ be any
positive (finite) measure on $X$, and assume that $m \times m$ restricted to
$X \times X\smallsetminus\Delta$ has as support all of $X \times
X\smallsetminus\Delta$.  Then one can consider all of the $L^p$-norms for $m
\times m$.  If one wants to put $\g$ (or $\rho$) explicitly into the picture,
one can consider the measure $\g(m \times m)$, although this just represents
the choice of a different measure.  Note that if $f$ is an ordinary Lipschitz
function for $\rho$, then $\g df$ is a bounded function on $X \times
X\smallsetminus\Delta$, so that $\|\g df\|_{p,m \times m}$ is finite.  Thus the
subalgebra of elements of $\cA$ for which this Lipschitz seminorm is finite is
dense in $\cA$.

To explore further possibilities, let us for simplicity assume that $X$ is
finite.  Then $\Omega_1^u = C(X \times X)$ can be viewed as the algebra of all
matrices whose entries are indexed by elements of $X \times X$.  The left and
right actions of $\cA$ on $\Omega_1^u$ can be viewed as coming from embedding
$\cA$ as the diagonal matrices and using left and right matrix
multiplication.  Then $\omega_0$ is the matrix with a $1$ in each entry.
On $\cA$ we keep the supremum norm, but on the matrix algebra $\Omega_1^u$ we
can consider any $\cA$-$\cA$-bimodule norm.  Let ${\Cal B}$ denote $\Omega_1^u$ viewed
as matrix algebra, and equipped with the usual $C^*$-algebra norm.  View
$\Omega_1^u$ as a ${\Cal B}$-${\Cal B}$-bimodule in the evident way.  
Then we can consider
${\Cal B}$-${\Cal B}$-bimodule norms on $\Omega_1^u$.  
Any such will in particular be an
$\cA$-$\cA$-bimodule norm.  But there has been extensive study of the possible
${\Cal B}$-${\Cal B}$-bimodule norms on $\Omega_1^u$.  
They are commonly called ``symmetric
norms'', and among the best known are the Schatten $p$-norms, which include
the Hilbert--Schmidt norm and the trace norm.  These have, of course, also
been extensively studied for operators on infinite dimensional Hilbert spaces,
and play a fundamental role in Connes' theory of integration on
non-commutative spaces. (See \cite{C2} Chapter~IV and its Appendix~D.  A nice
treatment of the finite case can be found in \cite{Bh}.)  From every symmetric
norm we obtain a Lip-norm on $\cA$ (since $\cA$ is finite-dimensional).  This
does not exhaust the possibilities, as there is no necessity to restrict to
symmetric norms in order to get $\cA$-$\cA$-bimodule norms.

All of the above discussion has been for the universal differential calculus.
We get many more possibilities by using other differential calculi.  We
continue to concentrate on the case of $\cA = C(X)$ with $X$ compact.  Now
sub-$\cA$-$\cA$-bimodules of $C(X\times X)$, when closed in the supremum norm,
will be ideals of $C(X \times X)$, and the quotient can be identified with
$C(W)$ for some closed subset $W$ of $X \times X$.  We can restrict $df$ to
$W$.  But some condition must be placed on $W$ if we want to ensure that
$df|_W = 0$ only if $f$ is a constant function.  For this purpose it is
convenient to assume, to begin with, that $W$ contains the diagonal $\Delta$
and is symmetric about $\Delta$, that is, if $(x,y) \in W$ then $(y,x)\in W$.
Given $x \in X$ we define the $W$-neighborhood of $x$ to be the (closed) set
of those $y \in X$ such that $(x,y) \in W$.  By the $W$-component of $x$ we
mean the smallest closed subset of $X$ which contains the $W$-neighborhood of
each of its points.  If $df|_W = 0$, then $f$ is constant on the $W$-component
of each point.  Thus a sufficient condition under which $df|_W = 0$ will imply
that $f$ is constant, is that the $W$-component of each point is all of $X$.
If $X$ is a finite set, then $W\smallsetminus\Delta$ can be viewed as consisting
of the directed edges for a graph whose vertices are the points of $X$.  Then
the above condition becomes the condition that this graph is connected in
the usual sense.  If $X$ is not discrete, it is usual to require that $W$ is
a neighborhood of $\Delta$.  Then each $W$-neighborhood of a point will be
an ordinary (closed) neighborhood, and so the $W$-component of each point will
be both closed and open.  In particular, if $X$ is connected it will be true
that $df|_W = 0$ implies that $f$ is constant.

We remark that if $W$ is a neighborhood of $\Delta$ and is symmetric about
$\Delta$, and if we set $\Omega = C(W)$, then the first order calculus
$(\Omega,d)$ obtained as above is the typical degree-one piece of the
complexes $(\Omega_W^*,d)$ used in defining the Alexander--Spanier cohomology
of $X$.  The higher-degree pieces are defined similarly but in terms of $X^n$
for various $n$.  The Alexander--Spanier cohomology is then obtained by taking
a limit of the homology of these complexes 
as $W$ shrinks to $\Delta$.  Essentially this view
can be seen in lemma $1.1$ of \cite{CM}, where smooth functions on a manifold
are used, and in Section~1 of \cite{MW}, where continuous functions are used.

Suppose now that $\Omega = C(W)$ as above, but assume now for simplicity that
$W$ and $\Delta$ are disjoint (with $W$ no longer required closed).  
Let $d$ be defined by $df = df|_W$, and assume
that if $df = 0$ then $f$ is a scalar multiple of $1$.  To obtain a Lipschitz
seminorm on $\cA$ we again just need to put a bimodule norm on $\Omega$.  The
method which is closest to the usual Lipschitz norm is to specify a nowhere
zero
function $\g$ on $W$ and set
$$
L(f) = \|\g df\|_{\infty}
$$
(on $W$, allowing value $+\infty$).  In this context however, if we set $\rho =
\g^{-1}$, it no longer makes much sense to ask that the triangle inequality
hold for $\rho$.  About the most that is reasonable is to ask that $\rho$,
hence $\g$, be positive, and that $\g(x,y) = \g(y,x)$ for $(x,y) \in W$, $x
\ne y$.  This is a situation which has been widely studied.  Entire books
\cite{Ra, RR} have been written about the problem of finding the corresponding
distance between two probability measure on $X$, often under the heading of
``the mass transportation problem''.  The function $\rho$ is then often called
a ``cost function''.  We should clarify that when $\rho$ is not a metric we
are dealing here with mass transportation ``with transshipment permitted''
\cite{RR}, not the original Monge--Kantorovich \cite{KA} mass
transportation problem, which does not permit transshipment, and may well not
yield a metric.  When transshipment is permitted and $\rho$ is not a metric on
$X$, the corresponding metric on $S(X)$ is called the Kantorovich--Rubenstein
metric \cite{KR1, KR2}.  
For a fascinating survey of some recent developments concerning the
original Monge--Kantorovich problem see \cite{Ev}.

When $X$ is a finite set and $W$ is viewed as specifying edges for a graph
which has $X$ as set of vertices, the cost function $\rho$ is naturally
interpreted as assigning lengths to the edges (though we will see a quite
different interpretation in Section 12).  Then the metric on $X$ coming
from $L_{\rho}$ is the usual path-length distance on the graph.  There has
been much study of how to compute this path-length distance efficiently for
large graphs.  We remark that if one prefers to have $\rho$ defined on all of
$X \times X$ one can simply set it equal to $+\infty$ on any $(x,y)$, $x \ne
y$, which is not an edge.

We remark that in the context of cost functions on compact sets there may well
be no non-constant functions for which the Lipschitz seminorm is finite.  As
one example let $X$ be the unit interval $[0,1]$, and set $\rho(x,y) =
|x-y|^2$.  This is, in effect, because we permit transshipment --- the original
Monge--Kantorovich problem is quite interesting for this particular cost
function, as shown in \cite{Ev}.  It is just that the minimal cost of moving
one probability measure directly to another does not then give a metric on
probability measures, because it may be less costly to use two or more moves.

There is a variety of other bimodule norms, such as
$L^p$-norms, which one can use 
for various differential calculi, and these give a wide variety
of metrics on probability measures \cite{Ra}.  A particularly deep application of such
norms, for the case of graphs, and involving explicitly Connes ideas of
non-commutative metrics, appears in \cite{Da}.  (I thank Nik Weaver for
bringing this paper to my attention.)

Let us now discuss briefly the case in which we have $\cA = M_n$, a full
matrix algebra.  As mentioned much earlier, one natural Lip-norm on $\cA$ is
just $L = \|\quad\|^{\sim}$.  Now $\cA'$ can be identified by means of the
normalized trace, $\tau$, with $\cA$ itself, but equipped with the
trace-norm.  Then ${\cA'}^0$, as in our earlier notation, consists of the
matrices with trace $0$.  Of course, $S(\cA)$ is identified with the positive
matrices of normalized trace $1$.  With this identification, we have
$$
\rho_L(\mu,\nu) = \text{trace}(|\mu-\nu|).
$$
This is exactly one of the metrics listed (with references) in the
introduction to \cite{ZS}.  Another one listed there uses the Hilbert--Schmidt
norm instead of the trace norm.  Listed also is a variety of other metrics on
$S(M_n)$ which have appeared in various applications.  But I have not checked
whether they come from Lip-norms.  There has also been much study of the
differential geometry of $S(M_n)$ for a variety of Riemannian metrics, 
especially the ``monotone metrics'', which are closely related
to operator monotone functions.  Two
very recent articles which contain many references to previous work on this
topic are \cite{Di, S}.  But the emphasis of most of this work is not on the
ordinary metric which a Riemannian metric induces on $S(M_n)$, but
rather on the differential geometric aspects. There
is also study 
of the volume form which is induced, and on associated probabilistic aspects.
For recent related study going in the direction of non-commutative entropy see
\cite{LR}.

\bigskip
\noindent
{\bf 11.  Dirac operators and differential calculi}

\medskip
We continue our comments of the previous section, but here we focus on how
Dirac operators fit into the picture.  Let $\cA$ be a unital $*$-algebra
equipped with a $C^*$-norm (perhaps not complete), and let $\pi$ be a faithful
representation of $\cA$, that is, an isometric $*$-homomorphism of $\cA$ into
the algebra $B(\cH)$ of bounded operators on a Hilbert space $\cH$.  Let $D$
be an essentially self-adjoint, possibly unbounded, operator on $\cH$, and
assume that $\pi(a)$ carries the domain of $D$ into itself for each $a \in
\cA$, and that on this domain $[D,\pi(a)]$ is a bounded operator, and so
extends uniquely to a bounded operator on $\cH$.  Then, following Connes, we
set
$$
L(a) = \|[D,\pi(a)]\|.
$$
As we did earlier, it is natural to require that $[D,\pi(a)] = 0$ only when
$a$ is a scalar multiple of $1$.  Many important examples of this situation
are now known.  But in general it seems difficult to ascertain whether the
corresponding metric on states gives the weak-$*$ topology, though this has
been shown for certain examples in \cite{Rf}. See also \cite{W2, W3, W5},
where the sets ${\Cal B}_t$ defined at the beginning of Section 3 are
shown to be totally bounded, in fact compact, 
for various examples.  We do not deal with
this question here, but rather try to relate the bimodule picture to the Dirac
picture.  One direction is apparent.  We view $B(\cH)$ as an
$\cA$-$\cA$-bimodule by setting
$$
aTb = \pi(a)T\pi(b).
$$
Then, although $D$ is only affiliated with $B(\cH)$, conceptually we use the
inner derivation which $D$ defines, so that
$$
da = D\pi(a) - \pi(a)D = [D,\pi(a)].
$$
(This, of course, is the starting point for Connes' non-commutative
differential calculus \cite{C2}.)  We then note that the operator norm on
$B(\cH)$ is an $\cA$-$\cA$-bimodule norm, and so upon setting
$$
L(a) = \|[D,\pi(a)]\|
$$
we obtain a Lipschitz norm, which we showed to be lower semicontinuous in
Proposition $3.8$.

But suppose we are given instead some first order differential calculus
$(\Omega,d)$ and a bimodule norm on $\Omega$ so that we obtain the
corresponding Lipschitz norm $L$.  Can we also obtain $L$ from a Dirac
operator?  For this to be possible we must have $L(a^*) = L(a)$, and $L$ must
be lower semicontinuous.  As mentioned earlier, $L$ must also fit into a family 
of ``matrix Lipschitz seminorms''. These conditions are probably not enough in
general, though I have not tried to find a counterexample.  But the following
superficial comments help to give some perspective.  (In most of the
considerations which follow the algebra structure on $\cA$ is only used in
order to get the Leibniz inequality.  Thus much of what follows actually works
for order-unit spaces.)

We saw in Proposition $10.3$ that we can extend $(\Omega,d)$ to obtain an
inner first-order calculus.  In analogy with this idea, suppose that we can
realize $\Omega$ as a subspace of $B(\cH)$ for some Hilbert space $\cH$, in
such a way that the norm on $\Omega$ is the operator norm, and the bimodule
structure is given by two $*$-representations, $\pi_1$ and $\pi_2$, of $\cA$
on $\cH$, so that
$$
a\omega b = \pi_1(a)\omega \pi_2(b)
$$
for $a,b \in \cA$ and $\omega \in \Omega$.  Suppose further that there is a
possibly-unbounded essentially self-adjoint operator, $D_0$, on $\cH$, such
that $\pi_1(a)$ and $\pi_2(a)$ carry the domain of $D_0$ into itself, and such
that
$$
da = D_0\pi_2(a) - \pi_1(a)D_0,
$$
which in particular must be a bounded operator.  Set $L(a) = \|da\|$.  This is
not exactly the Dirac operator setting, but it is not difficult to convert it
into that setting.  To arrange matters so that we have only one
representation, we let $\pi = \pi_1 \oplus \pi_2$ on $\cH \oplus \cH$ and set
$$
D_1 = \pmatrix 0 & D_0 \\ 0 & 0 \endpmatrix .
$$
Then we find that
$$
L(a) = \|[D_1,\pi(a)]\|.
$$
But of course $D_1$ is not self-adjoint.  We fix this in the traditional way
by again doubling the Hilbert space, with representation $\pi \oplus \pi$ of
$\cA$, and setting
$$
D = \pmatrix 0 & D_1^* \\ D_1 & 0 \endpmatrix .
$$
The corresponding Lipschitz norm is $L(a)\vee L(a^*)$, but from the
self-adjointness of $D$ one can check that we actually get back $L$.

Anyway, we are left with

\medskip
\noindent
{\bf 11.1\ Question}.  For an order-unit space $\cA$, or a $*$-algebra $\cA$
with $C^*$-norm, how does one characterize those Lip-norms on $\cA$ which come
from the Dirac operator construction?

\medskip
Even for finite-dimensional commutative $C^*$-algebras it is not clear to me
what the answer is.

As mentioned earlier, a Dirac operator also gives seminorms on all of the
matrix algebras over $\cA$, so that one can speak of this family as a ``matrix
Lipschitz norm'', in the spirit of \cite{Ef}.  Thus a related problem is to
characterize these structures.

Of course a given metric on $S(\cA)$ may come from several fairly different
Dirac operators.  For example, suppose that we have a compact space $X$, and a
closed neighborhood $W$ of the diagonal $\Delta$ of $X \times X$, together
with a cost function $\rho$ on $W$, just as in the previous section.  As
discussed there, we can use $\rho$ together with the first-order calculus
determined by $W$ to define a Lipschitz norm on $C(X)$.  (Further hypotheses
are needed for it to be a Lip-norm on a dense subalgebra of $C(X)$.)  Then by
the procedure discussed earlier in the present section we can pass to a Dirac
operator.  But that procedure enlarged the Hilbert space because a first-order
differential calculus usually involves two representations rather than one.
We will now show that there is an alternative method which does not enlarge the
Hilbert space.  This is a mild generalization of my lecture comments for
metric spaces mentioned earlier, whose details are indicated on page~274 of
\cite{W2}.  As earlier, let $m$ be a measure on $X$ of full support, and
consider $m \times m$ on $W\smallsetminus\Delta$.  Form the Hilbert space $\cH =
L^2(W\smallsetminus\Delta,m \times m)$.  We consider only the representation
$\pi$ of $\cA = C(X)$ on $\cH$ defined by
$$
(\pi_f\xi)(x,y) = f(x)\xi(x,y).
$$
(This is, of course, essentially the left action on the bimodule for $W$.)
Define an operator, $F$, on $\cH$ by the flip
$$
(F\xi)(x,y) = \xi(y,x).
$$
Because we are using a product measure, the operator $F$ is self-adjoint and
unitary.  Define an (unbounded) positive operator, $P$, on $\cH$ by
$$
(P\xi)(x,y) = \xi(x,y)/\rho(x,y).
$$
Because we assume that $\rho(x,y) = \rho(y,x)$ for all $(x,y) \in W$, the
operators $F$ and $P$ commute.  We define the Dirac operator by
$$
D = PF,
$$
so that $F$ is the phase of $D$ and $P = |D|$.  Informal calculation shows
that for any $f \in C(X)$ we have
$$
([D,\pi_f]\xi)(x,y) = ((f(y)-f(x))/\rho(x,y))\xi(y,x),
$$
so that
$$
L(f) = \|[D,\pi_f]\| = \sup\{|f(y)-f(x)|/\rho(x,y): (x,y) \in W\}.
$$
Of course, further hypotheses must be placed on $\rho$ in order for this to
give a Lip-norm.  But the right-hand side of the above equality is the usual
definition of a Lipschitz norm in this situation, especially in contexts such
as graph theory.  It will coincide with what one obtains in the corresponding
bimodule approach.  Notice that the resulting distance between two points
$x,y \in X$ can easily be strictly smaller than $\rho(x,y)$ (if $(x,y)$
happens to be in $W$).

For an interesting alternative (but closely related) method of obtaining the
usual distance on a graph (including infinite graphs) from a cost function, by
means of Dirac operators, see theorem $7.2$ of \cite{Da}.  Furthermore, in
\cite{Da} other very interesting and quite different Dirac operators
associated to cost functions on graphs are discussed in some detail, and used
to obtain improved estimates for heat kernels on graphs.  They can be
described in terms of first-order differential calculi and Laplace operators
along much the same lines as we used in Section~10.  Much of this is explicit
in \cite{Da}, and we will not elaborate on it here.

We should mention here that very interesting examples of Dirac operators
associated with non-commutative variants of sub-Riemannian manifolds appear in
the second example following axiom $4'$ of \cite{C3}, and in \cite{W5}.

\bigskip
\noindent
{\bf 12.  Resistance distance}

\medskip
We conclude with an appealing class of examples which do not fit into the
previous framework of differential calculi, and for which the Lip-norm does
not satisfy the Leibniz identity.  These examples come from graphs with ``cost
functions'' on the edges, but now the graph is interpreted as an electrical
circuit with resistances on the edges, whose values are given by the cost
function.  These examples have been extensively studied \cite{DS, Kl, KlR, KZ},
but I have not seen earlier mention of the corresponding metric on probability
measures which we will define here.  It is not clear to me whether this metric
is more than a curiosity.

All of the discussion here can be carried out for infinite graphs, along the
lines discussed extensively in \cite{DS}, but for simplicity we only discuss
finite graphs here.  The examples also have a fine alternative interpretation
in terms of random walks \cite{DS}.  Our term ``resistance distance'' is
taken from the title of \cite{KlR}.

The set-up, as indicated above, is a finite graph with set $X$ of vertices,
together with strictly positive real numbers $r_{xy} = r_{yx}$ assigned to
each (undirected) edge.  We interpret these numbers as resistances.  We assume
throughout that the graph is connected.  Given $x,y \in X$, $x \ne y$, we can
imagine putting a voltage difference across $x$ and $y$, adjusted so that one
unit of current flows in at $x$ and out at $y$.  Then Ohm's law says that the
``effective resistance'' is equal to the required voltage difference.  We
denote this effective resistance by $\rho(x,y)$.  It is, in fact, a metric on
$X$.  The only reference I know for this is \cite{KlR, K, KZ}, but my friends in
probability theory tell me that within the context of random walks rather than
resistances this is well-known, even if no reference comes to mind.

Suppose now that $\mu$ and $\nu$ are general probability measures on $X$.
Although it does not seem so intuitively obvious, we will see shortly that we
can establish voltages on the points of $X$ such that unit total current flows
into the circuit, with the amount flowing in at each point $x$ given by
$\mu_x$, while unit total current flows out of the circuit, with the amount at
each point given by $\nu$ (with the evident interpretation when the supports
of $\mu$ and $\nu$ are not disjoint).  For the analysis of this situation it
is useful to define a function, $c$, on the edges, by $c_{xy} = 1/r_{xy}$.
This is commonly called the ``conductance''.  It is convenient to extend $c$
to all of $X \times X$ by setting $c_{xy} = 0$ if $(x,y)$ is not an edge (or
if $y = x$).  Let $f \in C(X)$, interpreted as voltages applied to the points
of $X$.  We let $df$ be defined as earlier for the universal calculus (or for
the calculus corresponding to the edges).  We let $\nabla f$ denote the
resulting flow inside the circuit.  By Ohm's law the flow (before
electrons were discovered) from $x$ to $y$ is
given by
$$
(\nabla f)(x,y) = (f(x) - f(y))c_{xy} = -c(df),
$$
where by $c(df)$ we mean the pointwise product of functions.  Note that
$\nabla f$ is a function on directed edges, with
$$
(\nabla f)(x,y) = -(\nabla f)(y,x)
$$
(and value $0$ if $(x,y)$ is not an edge).

Suppose now that $\omega$ is any function on directed edges such that
$\omega(x,y) = -\omega(y,x)$.  We interpret $\omega(x,y)$ as giving the
magnitude of a current from $x$ to $y$.  (To be more realistic we should
require $0$ circulation, but we will have no need to impose this
requirement.)  To sustain this current, we will in general have to insert (or
extract) current at various vertices.  We let $\text{div}(\omega)(x)$ denote
the current which must be inserted at $x$.  By Kirchhoff's laws we have
$$
\text{div}(\omega)(x) = \sum_y \omega(x,y).
$$
Note that because $\omega(x,y) = -\omega(y,x)$, we will have
$$
\sum_x \text{div}(\omega)(x) = 0,
$$
which accords with the fact that the total amount of current inserted must be
$0$.

Suppose now that $f \in C(X)$ and that we set $\omega = \nabla f$.  We see
from above that the currents which must be inserted to sustain the voltages
given by $f$ must be
$$
\text{div}(\nabla f),
$$
which we denote by $\Delta f$.  To accord with our earlier notation, we let
${\cA'}^0$ denote the signed measures, $\lam$, on $X$ for which $\<1,\lam\> =
0$.  The discussion of the previous paragraph can be
interpreted as saying that $\Delta f \in
{\cA'}^0$.

Suppose now that we are given $\lam \in {\cA'}^0$.  Can we find $f$ such that
$\Delta f = \lam$?  Note that since $\Delta 1 = 0$, we know that $f$ will not
be unique, but rather that, as usual with potential functions, we can expect $f$ to be
unique only up to a constant function.  To proceed further we must more
carefully analyze the operator $\Delta$ in the traditional way \cite{DS, K}.
For $f \in C(X)$ we have
$$
\align
(\Delta f)(x) &= \sum_y (\nabla f)(x,y) \\
&= \sum_y (f(x)-f(y))c_{xy} 
= f(x) \sum_y c_{xy} - \sum_y f(y)c_{xy}.
\endalign
$$
Let $D$ denote the diagonal matrix with diagonal entries
$$
D_{xx} = \sum_y c_{xy}.
$$
If we view $f$ as a column vector, we see that
$$
\Delta f = (D-C)f.
$$

From the Peron--Frobenius theorem and the fact that our graph is connected, it
follows that the kernel of $\Delta$ consists exactly of the constant
functions. If we permit ourselves to confuse vector spaces a bit, we see
that $\Delta$ is self-adjoint with respect to the standard inner-product
on column vectors. Thus it carries the orthogonal complement, $\Cal H$,
of the constant functions into itself, and it is invertible on $\Cal H$.
Consequently, for every $\lam \in {\cA'}^0$ we can find a unique $f \in {\Cal H}$
such that $\Delta f = \lam$. We will write this as $f = \Delta^{-1}\lam$, where
we view $\Delta$ as restricted to $\Cal H$ so that it is invertible there.

Suppose now that $x$ and $y$ are fixed points of $X$, and that $\lam = \de_x -
\de_y$, where $\de_x$ denotes the $\de$-measure at $x$.  Thus we are inserting
one unit of current at $x$ and extracting it at $y$.  Let $f =
\Delta^{-1}\lam$.  According to our earlier comments, the effective resistance
from $x$ to $y$, $\rho(x,y)$, is given by $f(x)-f(y) = (\Delta^{-1}\lam)(x) -
(\Delta^{-1}\lam)(y)$.  It is now easy to see why $\rho$ is a metric, along
the lines given in \cite{KlR}.  If $z$ is any other point of $X$, let
$$
g = \Delta^{-1}(\de_x - \de_z),\qquad h = \Delta^{-1}(\de_z -\de_y).
$$
Clearly $f = g+h$, so
$$
\rho(x,y) = g(x) - g(y) + h(x) - h(y).
$$
But simple considerations show that $g$ must take its maximum and minimum
values at $x$ and $z$, so that
$$
g(x) - g(y) \le g(x) - g(z) = \rho(x,z).
$$
Similarly $h(x) - h(y) \le \rho(x,z)$.  The triangle inequality for $\rho$
follows.

But we are interested more generally in the effective resistance between $\mu$
and $\nu$ where $\mu$ and $\nu$ are arbitrary probability measures, and it is
not even clear how this should be defined.  (It does not seem natural just to
use the Monge--Kantorovich metric from $\rho$.)  In view of our earlier
considerations we should form $\lam= \mu-\nu$, and so we need an appropriate
norm on ${\cA'}^0$, and this should be the dual norm of a Lip-norm, say $L$,
on $C(X)$, probably defined by means of a norm on $\Omega^u$.  The dual
norm, $L'$, should be such that if $\lam = \de_x -\de_y$, then $L'(\lam) =
(\Delta^{-1}\lam)(x) - (\Delta^{-1}\lam)(y)$.  But as remarked above,
$\Delta^{-1}\lam$ takes its maximum and minimum values at $x$ and $y$.  Thus a
norm which will meet this requirement is
$$
L'(\lam) = 2\|\Delta^{-1}\lam\|_{\infty}^{\sim},
$$
where $\|\quad\|_{\infty}^{\sim}$ is as defined in Section~1.  To find $L$ on
$C(X)$ we use the self-adjointness of $\Delta$ to calculate, for $g \in C(X)$
and any $\lam \in {\cA'}^0$,
$$
\<g,\lam\> = \<g,\Delta\Delta^{-1}\lam\> = \<\Delta g,\Delta^{-1}\lam\>.
$$
The supremum over $\lam$ such that $2\|\Delta^{-1}\lam\|_{\infty}^{\sim} \le
1$ is the same as the supremum of
$$
\<(1/2)\Delta g,h\>
$$
over $h$ such that $\|{\tilde h}\|_{\infty}^{\sim} \le 1$.  But we saw earlier
that this gives just the restriction to ${\cA'}^0$ of the dual norm for
$\|\quad\|_{\infty}$ on $C(X)$, which is the $L^1$-norm.  Thus we see that we
must set
$$
\align
L(g) &= (1/2) \|\Delta g\|_1 = (1/2) \sum_x |(\Delta g)(x)| \\
&= (1/2) \sum_x \left| \sum_y (g(x)-g(y))c_{xy}\right|
= (1/2) \sum_x \left| \sum_y dg(x,y)c_{xy}\right|.
\endalign
$$
This is certainly rather different from the usual Lip-norms for metrics on
finite sets.  The above expression suggests that we define a seminorm, $N$, on
$\Omega^u$ by
$$
N(\omega) = (1/2) \sum_x \left| \sum_y \omega(x,y)c_{xy}\right|,
$$
so that we have
$$
L(g) = N(dg).
$$
Reversal of the earlier calculation shows that the dual norm is the $L'$
considered above, so that we obtain the desired $\rho(\mu,\nu)$.  However $N$
will not usually be a bimodule norm, so that we are not fully in the context
of the previous sections, and $L$ need not satisfy the Leibniz inequality.

I must admit that I see no particularly natural interpretation for $L(g)$, nor
for $\rho(\mu,\nu)$, even if we call the latter ``effective resistance''.  If
$g$ were interpreted as giving voltages on $X$, then $L(g)$ would be half the
sum of the absolute values of the currents inserted or extracted from the
circuit, and thus exactly the sum of the currents inserted into the circuit
(disregarding the currents extracted).  But I do not see why it is natural to
give $g$ such an interpretation as voltages.  If one goes back to the
effective resistance between two points, then it is easily seen that this is
equal to the energy dissipated by the circuit when one unit of current is
inserted.  This suggests using the dissipated energy in the more general case
of arbitrary probability measures $\mu$ and $\nu$.  But the energy dissipated
along any edge varies as the square of the current, and one can see by
examples that this causes the triangle inequality to fail.  One does obtain a
metric if one uses the square-root of the dissipated energy, but this does not
give the correct value for the effective resistance between two points.  These
possibilities are not far from the Lipschitz norm used right after lemma $4.1$
of \cite{Da} to define the metric denoted there by $d_3$.  This Lipschitz norm
can be interpreted as the supremum over the points $x$ of $X$ of the square
roots of the energy dissipations in all the edges beginning at $x$.  Perhaps
the discussion of Dirichlet spaces given in section~6 of \cite{W6}, or the
``twisted bimodule structure'' and corresponding differential discussed
beginning on page~149 of \cite{Me} in connection with Hudson's treatment of
discrete flows and stochastic differential equations, could be used to shed
more light on this.  Or perhaps some of the stopping rules or mixing times
considered for Markov chains, as discussed in \cite{LW}, are relevant.

Finally, we remark that it would be interesting to study resistance distance
in the continuous case, for example for thin plates of resistance metal of
various shapes.


\newpage
\Refs
\widestnumber\key{BCD}

\ref \key Al
\by Alfsen, E. M.
\book Compact convex sets and boundary integrals
\publ Springer-Verlag
\publaddr Berlin, New York
\yr 1971
\endref

\ref \key Ar
\by Arveson, W. B.
\paper The harmonic analysis of automorphism groups
\inbook Proc. Sympos. Pure Math.
\publ Amer. Math. Soc., Providence, RI \vol 38 \pages
199--269  \yr 1982
\paperinfo MR 84m:46085
\endref

\ref \key Ba
\by Baaj, S.
\paper  Multiplicateurs non born\'es
\paperinfo Th\`ese 3eme cycle, Universit\'e Paris VI \yr 1980
\endref

\ref \key BCD
\by Bade, W. G., Curtis, P. C., and Dales, H. G.
\paper Amenability and weak amenability for Beurling and Lipschitz
algebras
\jour Proc. London Math. Soc. (3) {\bf 55}, no.~2 \yr 1987 \pages
359--377
\paperinfo MR  88f:46098
\endref

\ref \key Bh
\by Bhatia, R.
\book Matrix Analysis
\publ Springer-Verlag
\publaddr New York, New York \yr 1997  \bookinfo MR
98i:15003
\endref

\ref \key BC
\by Blackadar, B. and Cuntz, J.
\paper  Differential Banach algebra norms and smooth
subalgebras of $C^*$-algebras
\jour   J. Operator theory      \vol    26
\yr     1991    \pages  255--282
\endref

\ref \key C1
\by Connes, A.
\paper Compact metric spaces, Fredholm modules and hyperfiniteness
\jour   Ergodic Theory and Dynamical Systems    \vol    9
\yr     1989    \pages  207--220
\endref

\ref \key C2
\by Connes, A.
\book   Noncommutative Geometry
\publ   Academic Press  \yr     1994
\publaddr       San Diego
\endref

\ref \key C3
\by Connes, A.
\paper Gravity coupled with matter and the foundation of
non commutative geometry
\jour Comm. Math. Phys. \vol    182
\yr     1996    \pages  155--176
\paperinfo  hep-th/9603053.
\endref

\ref \key CM
\by Connes, A. and Moscovici, H.
\paper  Cyclic cohomology, the Novikov conjecture and hyperbolic groups
\jour Topology {\bf 29}, no.~3
\yr 1990 \pages 345--388 \paperinfo MR 92a:58137
\endref

\ref    \key  Cw
\book  A Course in Functional Analysis, 2nd edition
\by Conway, J. B.
\publ Springer Verlag      \yr 1990
\publaddr    New York
\endref

\ref \key Da
\by Davies, E. B.
\paper Analysis on graphs and noncommutative geometry
\jour J. Funct. Anal. {\bf
111}, no.~2
\yr 1993
\pages 398--430 \paperinfo MR 93m:58110
\endref

\ref \key  Di
\by Dittmann, J.
\paper On the curvature of monotone metrics and a conjecture concerning
the Kubo-Mori metric
\paperinfo  quant-ph/9906009
\endref

\ref \key DS
\by Doyle, P. G. and Snell, J. L.
\book Random Walks and Electric Networks
\publ Mathematical Association of
America
\publaddr Washington, DC
\yr 1984 \bookinfo MR 89a:94023
\endref

\ref \key Ef
\by Effros, E. G.
\paper Advances in quantized functional analysis
\inbook Proceedings of the
International Congress of
Mathematicians, Vol. 1, 2 (Berkeley, Calif., 1986)
\publ Amer. Math.
Soc.,  Providence, RI \pages 906--916 \yr
1987 \paperinfo MR 89e:46064
\endref

\ref \key Ev
\by Evans, L. C.
\paper Partial differential equations and Monge--Kantorovich mass transfer
\jour preprint
\paperinfo  available at: http://berkeley.math.edu/~evans/
\endref

\ref \key H
\by Hanin, L. G.
\paper  Kantorovich-Rubinstein norm and its application in the theory of
Lipschitz spaces
\jour Proc. Amer.
Math. Soc. {\bf 115}, no.~2  \yr 1992 \pages 345--352 \paperinfo MR 92i:46026
\endref

\ref \key J
\by  Johnson, B. E.
\paper Amenability and weak amenability for Beurling and Lipschitz
algebras
\jour Proc. London Math. Soc. (3) {\bf 55}, no.~2 \yr 1987 \pages 359--377
\paperinfo
MR  88f:46098
\endref

\ref \key K1
\by Kadison, R. V.
\paper A representation theory for commutative topological algebra
\jour Mem. Amer.
Math. Soc. {\bf 1951}, no.~7 \yr 1951 \pages 39 pp. \paperinfo MR {\bf 13},
360b
\endref

\ref \key K2
\by Kadison, R. V.
\paper Unitary invariants for representations of operator algebras
\jour Ann. of
Math. (2) \vol 66 \yr 1957 \pages
304--379 \paperinfo MR {\bf 19}, 665e
\endref

\ref \key K3
\by Kadison, R. V.
\paper Transformations of states in operator theory and dynamics
\jour Topology {\bf
3}, suppl.~2 \yr 1965 \pages
177--198 \paperinfo MR {\bf 29} \#6328
\endref

\ref \key KA
\by Kantorovi\v c, L. V. and Akilov, G. P.
\book  Functional Analysis (Russian)
\publ Third edition, ``Nauka'' \publaddr Moscow \yr 1984 \bookinfo MR
86m:46001
\endref

\ref \key KR1
\by Kantorovi\v c L. V. and Rubin\v ste\u\i n, G. \v S.
\paper On a functional space and certain extremum problems
\jour Dokl.
Akad. Nauk SSSR (N.S.) \vol 115 \yr 1957 \pages 1058--1061 \paperinfo MR {\bf
20}
\#1219
\endref

\ref \key KR2
\by Kantorovi\v c L. V. and Rubin\v ste\u\i n, G. \v S.
\paper On a space of completely
additive functions
\jour Vestnik Leningrad. Univ. {\bf 13}, no.~7  \yr 1958 \pages 52--59
\paperinfo MR {\bf 21} \#808
\endref

\ref \key Kl
\by Klein, D.J.
\paper Graph geometry, graph metrics, and Wiener
\jour Match No.~35 \yr 1997 \pages 7--27 \paperinfo
MR  98e:05103
\endref

\ref \key KlR
\by Klein, D. J. and Randi\'c, M.
\paper Resistance distance
\jour J. Math. Chem. {\bf 12}, no.~1-4 \yr 1993 \pages 81--95 \paperinfo
MR  94d:94041
\endref

\ref \key KZ
\by Klein, D. J. and Zhu, H.-Y.
\paper Distances and volumina for graphs
\jour J. Math. Chem. {\bf 23}, no.~1-2 \yr 1998 \pages 179--195 \paperinfo MR
99f:05032
\endref

\ref \key LR
\by Lesniewski, A. and Ruskai, M. B.
\paper Monotone Riemannian metrics and
relative entropy on non-commutative probability spaces
\paperinfo math-ph/9808016.
\endref

\ref \key LW
\by Lov\'asz L. and Winkler, P.
\paper Mixing times
\inbook Microsurveys in Discrete Probability (Princeton,
NJ)
\publ Amer. Math. Soc., Providence, RI
\yr 1997 \pages 85--133
\endref

\ref \key Me
\by  Meyer, P.-A.
\paper Quantum probability for probabilists
\jour Lecture Notes in Math., 1538,
Springer-Verlag, Berlin \yr 1993 \paperinfo MR 94k:81152
\endref

\ref \key MW
\by Moscovici, H. and Wu, F.-B.
\paper Index theory without symbols
\inbook $C^*$-Algebras:  1943--1993
(San Antonio,
TX, 1993)
\jour Contemp. Math. \vol 167
\publ Amer. Math. Soc., Providence,
RI \pages 304--351 \yr 1994 \paperinfo MR
96g:58184
\endref

\ref \key P
\by Pavlovi\'c, B.
\paper Defining metric spaces via operators from unital $C^*$-algebras
\jour
Pacific J.  Math. {\bf 186}, no.~2 \yr 1998 \pages 285--313
\endref

\ref \key Ra
\by Rachev, S. T.
\book Probability Metrics and the Stability of Stochastic Models
\publ John
Wiley and  Sons \yr 1991
\endref

\ref \key RR
\by Rachev, S. T. and R\"uschendorf, L.
\book Mass Transportation Problems. Vol.
I, Theory
\publ Springer-Verlag \publaddr New York,
New York \yr 1998
\endref

\ref \key Rf
\by Rieffel, M. A.
\paper Metrics on states from actions of compact groups
\jour  Doc. Math. \vol 3 \yr 1998 \pages 215--229
\paperinfo math.OA/9807084
\endref

\ref \key Ru
\by Rudin, W.
\book Functional Analysis
\publ Second edition, McGraw-Hill, Inc. \publaddr New York, New York \yr
1991  \paperinfo MR 92k:46001
\endref

\ref \key S
\by Slater, P. B.
\paper  A priori probabilities --- based on volume elements of monotone
metrics  --- of quantum disentanglements
\paperinfo quant-ph/9810026
\endref

\ref \key W1
\by Weaver, N.
\book Lipschitz Algebras
\publ  World Scientific    \yr 1999
\publaddr   Singapore
\endref

\ref \key W2
\by Weaver, N.
\paper Lipschitz algebras and derivations of von Neumann algebras
\jour J. Funct.
Anal. {\bf 139}, no.~2 \yr 1996 \pages 261--300 \paperinfo MR 97f:46081
\endref

\ref \key W3
\by Weaver, N.
\paper $\alpha$-Lipschitz algebras on the noncommutative torus
\jour J. Operator
Theory \vol 39 \yr 1998 \pages
123--138
\endref

\ref \key W4
\by Weaver, N.
\paper  Operator spaces and noncommutative metrics
\jour      preprint
\endref

\ref \key W5
\by Weaver, N.
\paper Sub-Riemannian metrics for quantum Heisenberg
manifolds
\paperinfo      math.OA/9801014
\endref

\ref \key W6
\by Weaver, N.
\paper Lipschitz algebras and derivations, II:  Exterior differentiation
\paperinfo math.FA/9807096
\endref

\ref \key Wo
\by Woronowicz, S. L.
\paper Unbounded elements affiliated with $C^*$-algebras and noncompact
quantum groups
\jour Comm. Math. Phys. {\bf 136}, no.~2 \yr 1991 \pages 399--432
\paperinfo MR 92b:46117
\endref

\ref \key ZS
\by Zyczkowski, K. and W. S\l omczy\'nski, W.
\paper The Monge distance between quantum states
\jour J. Phys. A {\bf 31}, no.~45 \yr 1998 \pages 9095--9104
\paperinfo quant-ph/9711011
\endref

\endRefs
\enddocument